\documentclass[letterpaper,11pt]{amsart}

\usepackage[utf8]{inputenc}
\usepackage{amsmath,amssymb,amsthm,mathrsfs}

\usepackage{array}
\usepackage{tikz-cd}
\usepackage[colorlinks=true,linkcolor=blue,citecolor=blue]{hyperref}
\usepackage[shortlabels]{enumitem}

\setlist[enumerate]{
  label=(\alph*),
  itemsep=0.25em
}
\numberwithin{equation}{section}

\let\emptyset\varnothing

\usepackage{marginnote} 
\usepackage{booktabs}

\theoremstyle{plain}
\newtheorem{theorem}{Theorem}[section]
\newtheorem{proposition}[theorem]{Proposition}
\newtheorem{lemma}[theorem]{Lemma}
\newtheorem{corollary}[theorem]{Corollary}

\theoremstyle{plain}
\newtheorem{maintheorem}{Theorem}

\theoremstyle{definition}
\newtheorem{definition}[theorem]{Definition}
\newtheorem{example}[theorem]{Example}

\theoremstyle{remark}
\newtheorem{remark}[theorem]{Remark}

\newcommand{\K}{\mathbb{K}}
\newcommand{\T}{\mathbb{T}}
\newcommand{\R}{\mathbb{R}}
\newcommand{\C}{\mathbb{C}}
\newcommand{\N}{\mathbb{N}}
\newcommand{\Lip}{\operatorname{Lip}}
\newcommand{\LipNRA}{\operatorname{LipNRA}}

\newcommand{\Lipz}{\operatorname{Lip}_0}

\newcommand{\Lin}{\mathcal{L}}
\newcommand{\F}{\mathcal{F}}
\newcommand{\id}{\operatorname{Id}}
\newcommand{\Rea}{\operatorname{Re}}
\newcommand{\vL}[1]{\nu_{L,#1}}

\DeclareMathOperator{\spann}{span}

\DeclareMathOperator{\conv}{conv}
\newcommand{\restricted}{\mathord{\upharpoonright}}

\begin{document}

\title[The numerical index with respect to a Lipschitz map]
{The numerical index with respect to a Lipschitz map} 

\author[G.~Choi]{Geunsu Choi}
\address[G.~Choi]{Department of Mathematics Education, Sunchon National University, 57922 Suncheon, Republic of Korea\newline
\href{http://orcid.org/0000-0002-4321-1524}{ORCID: \texttt{0000-0002-4321-1524}  }}
\email{\texttt{gschoi@scnu.ac.kr}}

\author[M.~Han]{Manwook Han}
\address[M.~Han]{G-LAMP Research Center, Chungbuk National University, 28644 Cheongju, Republic of Korea\newline
\href{https://orcid.org/0009-0000-7250-1814}{ORCID: \texttt{0009-0000-7250-1814}  }}
\email{\texttt{626086@chungbuk.ac.kr}}

\author[M.~Jung]{Mingu Jung}
\address[M.~Jung]{Department of Mathematics \& Research Institute for Natural Sciences, Hanyang University, 04763 Seoul, Republic of Korea \newline 
\href{http://orcid.org/0000-0003-2240-2855}{ORCID: \texttt{0000-0003-2240-2855} } }
\email{\texttt{mingujung@hanyang.ac.kr}}

\keywords{Lipschitz numerical index, Lipschitz numerical radius, numerical index with respect to an operator, Lipschitz-free spaces, graph metrics, bipartite graphs, numerical-radius attainment}
\subjclass[2020]{Primary 47A12; Secondary 46B04, 46B20, 05C12}

\date{\today}

\begin{abstract}
We introduce the Lipschitz numerical range, radius, and index with respect to a norm-one Lipschitz map, unifying the Lipschitz numerical index of a Banach space and the numerical index with respect to an operator. For real-valued maps, we obtain metric formulas on finite metric spaces and graphs. In particular, a connected simple unweighted graph with at least two vertices is bipartite if and only if its real Lipschitz-free space, or equivalently its dual, has classical numerical index one. We also characterize finite metric spaces by universal numerical-radius attainment. Among complete metric spaces, we characterize length spaces by non-denseness of maps attaining their $g$-Lipschitz numerical radius for every norm-one $g$. Moreover, we derive upper bounds for the $g$-Lipschitz numerical index with $\ell_p$ targets, $1<p<\infty$, and lower bounds for the Lipschitz numerical index of spaces with property $\beta$. Finally, for norm-one operators $G: X \to Y$ between real Banach spaces, we prove
$$
n_{G^*}(Y^*,X^*) \leq n_{L,G}(X,Y) \leq n_G(X,Y),
$$
with equality in the upper bound whenever $X$ is separable or $Y$ is $1$-complemented in its bidual. This extends recent equality results for the identity operator to arbitrary norm-one operators and a broader class of target spaces.
\end{abstract}

\maketitle

\tableofcontents

\section{Introduction}\label{sec:intro}

\subsection{Background and motivation}

Throughout, $\K$ denotes either $\R$ or $\C$, and
$\T=\{\theta\in\K:|\theta|=1\}$. For a Banach space $X$, we write
$B_X$, $S_X$, and $X^*$ for its closed unit ball, unit sphere, and
dual space, respectively. We denote by $J_X\colon X\to X^{**}$ the canonical isometric embedding. For $\lambda\ge1$, a closed subspace $Z$ of $X$ is said to be \emph{$\lambda$-complemented} in $X$ if there exists a bounded linear projection $P\colon X\to X$ with $P(X)=Z$ and $\|P\|\le\lambda$. We denote by $\Lin(X,Y)$ the space of bounded linear operators from $X$ to $Y$, and write
$\Lin(X)=\Lin(X,X)$.

The numerical index of a Banach space measures the extent to which the
numerical radius is equivalent to the operator norm. Given a Banach
space $X$, let
\[
\Pi(X)
=
\left\{
(x,x^*)\in S_X\times S_{X^*}:x^*(x)=1
\right\}.
\]
The numerical radius of $T\in\Lin(X)$ and the numerical index of $X$
are defined by
\[
\nu(T)
=
\sup\left\{
|x^*(Tx)|:(x,x^*)\in\Pi(X)
\right\}
\]
and
\[
n(X)
=
\inf\left\{
\nu(T):T\in\Lin(X),\ \|T\|=1
\right\},
\]
respectively. We refer to the monograph \cite{BD} for the classical
theory of numerical ranges and numerical indices.

Two extensions of this notion motivate the present work. The first is
obtained by replacing the identity operator with an arbitrary
norm-one operator. Ardalani \cite{Ard} initiated the study of the
numerical radius and numerical index with respect to a norm-one
operator $G\in\Lin(X,Y)$, and the resulting theory was developed
systematically in \cite{KMMPQ}. For $0<\delta<1$, one considers
\[
\Pi_G^{\mathrm{lin}}(\delta)
=
\left\{
(x,y^*)\in S_X\times S_{Y^*}:
\Rea y^*(Gx)>1-\delta
\right\}
\]
and, for $T\in\Lin(X,Y)$, defines
\[
\nu_G^\delta(T)
=
\sup\left\{
|y^*(Tx)|:(x,y^*)\in\Pi_G^{\mathrm{lin}}(\delta)
\right\},
\quad
\nu_G(T)
=
\inf_{0<\delta<1}\nu_G^\delta(T).
\]
The numerical index with respect to $G$ is then
\[
n_G(X,Y)
=
\inf\left\{
\nu_G(T):T\in\Lin(X,Y),\ \|T\|=1
\right\}.
\]
For $X=Y$ and $G=\id_X$, this recovers the classical numerical index
$n(X)$. Moreover, $n_G(X,Y)=1$ if and only if $G$ is a spear
operator, that is,
\[
\max_{\theta\in\T}\|G+\theta T\|
=
1+\|T\|
\quad
(T\in\Lin(X,Y));
\]
see \cite{KMMP} for a systematic treatment of spear operators.

The second extension is obtained by replacing linear operators with
Lipschitz maps. Wang, Huang and Tan \cite{WHT} introduced the
Lipschitz numerical radius and the Lipschitz numerical index of a
Banach space. If $\Lipz(X,X)$ denotes the space of Lipschitz maps
$f:X\to X$ satisfying $f(0)=0$, then
\[
\nu_L(f)
=
\sup\left\{
\frac{|x^*(f(x)-f(y))|}{\|x-y\|^2}:
x\neq y,\ x^*\in D(x-y)
\right\},
\]
where
\[
D(x)
=
\left\{
x^*\in X^*:
x^*(x)=\|x^*\|\|x\|=\|x\|^2
\right\},
\]
and
\[
n_L(X)
=
\inf\left\{
\nu_L(f):f\in\Lipz(X,X),\ \|f\|_L=1
\right\}.
\]
Since $\nu_L$ agrees with the classical numerical radius on
$\Lin(X)$, one always has $n_L(X)\leq n(X)$.

Several positive equality results are known in the real case. If $X$
is a real lush space, then $n_L(X)=n(X)=1$
\cite[Theorem~2.6]{WHT}. More recently, P\'erez-Hern\'andez proved
that
\[
n_L(X)=n(X)
\]
whenever $X$ is a real separable Banach space or a real dual Banach
space; see \cite[Theorems~2.1 and~3.1]{PH}. Further results concerning
the range of values of $n_L$ under renorming and Lipschitz
numerical-radius attaining maps were obtained in \cite{CJT}.

The classical Lipschitz numerical radius tests the increments of a map
against functionals norming the corresponding increments of the identity. This suggests replacing the identity by a fixed norm-one Lipschitz map $g:M\to Y$, where $M$ is a pointed metric space and $Y$ is a Banach space. The resulting question is how much of the Lipschitz norm of another map can be recovered from the pairs and functionals on which $g$ almost attains its norm.

A connection between Lipschitz numerical ranges and numerical ranges
with respect to linear operators is already available through Lipschitz-free linearization; see \cite[Section~5]{KMMPQ}. We develop this viewpoint for arbitrary norm-one Lipschitz maps $g$ and investigate how the resulting numerical quantities reflect the metric structure of the domain and the geometry of the target space.

\subsection{The $g$-Lipschitz numerical radius and index}

Throughout this paper, $M$ denotes a pointed metric space with at least two points, and $0$ denotes its distinguished point. All Banach spaces are assumed to be nonzero and over a scalar field $\mathbb{K} \in \{\mathbb{R},\mathbb{C}\}$. Let $Y$ be a Banach space. We denote by $\Lipz(M,Y)$ the Banach space of all Lipschitz maps $f:M\to Y$ satisfying $f(0)=0$, endowed with the \emph{Lipschitz norm}
\[
\|f\|_L
=
\sup_{p\neq q}
\frac{\|f(p)-f(q)\|}{d(p,q)}.
\]

A norm-one Lipschitz map need not attain its Lipschitz norm at any pair of distinct points. We therefore use pairs and functionals on which its normalized increments are almost normed, and then pass to the limit as the error tends to zero.

\begin{definition}\label{def:main}
Let $M$ be a pointed metric space, $Y$ be a Banach space, and let
$g\in\Lipz(M,Y)$ satisfy $\|g\|_L=1$. For $0<\delta<1$, set
\[
\Pi_g(\delta) =
\left\{
(p,q,y^*)\in M\times M\times S_{Y^*}: p\neq q,  \Rea y^* \left( \frac{g(p)-g(q)}{d(p,q)} \right) >1-\delta
\right\}.
\]
For $f\in\Lipz(M,Y)$, define
\[
\vL{g}^\delta(f)
=
\sup\left\{
\left| y^* \left( \frac{f(p)-f(q)}{d(p,q)} \right) \right|: (p,q,y^*)\in\Pi_g(\delta)
\right\}
\]
and
\[
W_{L,g}^{\delta}(f)
= \left\{ y^* \left( \frac{f(p)-f(q)}{d(p,q)} \right) : (p,q,y^*)\in\Pi_g(\delta) \right\}.
\]
The \emph{$g$-Lipschitz numerical radius} of $f$ and the
\emph{$g$-Lipschitz numerical range} of $f$ are defined, respectively,
by
\[
\nu_{L,g}(f)
= \inf_{0<\delta<1}\vL{g}^\delta(f), \quad \text{and} \quad W_{L,g}(f)
= \bigcap_{0<\delta<1}
\overline{W_{L,g}^{\delta}(f)}.
\]
Finally, the \emph{$g$-Lipschitz numerical index} of the pair $(M,Y)$
is
\[
n_{L,g}(M,Y)
=
\inf\left\{ \nu_{L,g}(f):
f\in\Lipz(M,Y),\ \|f\|_L=1
\right\}.
\]
\end{definition}

Since $\|g\|_L=1$, the set $\Pi_g(\delta)$ is nonempty for every
$0<\delta<1$. Moreover,
\[
\Pi_g(\delta)\subseteq\Pi_g(\delta')
\quad
(0<\delta\leq\delta'<1),
\]
and hence $\nu_{L,g}(f) = \lim_{\delta\to0^+}\vL{g}^\delta(f).$ By homogeneity,
\[
n_{L,g}(M,Y)\|f\|_L \le \nu_{L,g}(f) \le \|f\|_L \quad(f\in\Lipz(M,Y)).
\]
Moreover, $n_{L,g}(M,Y)$ is the largest constant for which the left-hand inequality holds for every $f\in\Lipz(M,Y)$. Remark~\ref{rem:linear} relates the preceding definition to the two motivating theories.

\begin{remark}\label{rem:linear}
Let $X$ and $Y$ be Banach spaces.
\begin{enumerate}[(a)]
\item If $g=\id_X$, then
\vspace{0.25em}

\begin{itemize}
\itemsep0.25em
    \item $\nu_{L,g}(f)=\nu_L(f)$ for every $f \in \Lip_0 (X,X)$, and 
    \item $n_{L,g}(X,X)=n_L(X)$.
\end{itemize}
\vspace{0.25em}

\noindent These identities follow from the characterization of the classical Lipschitz numerical range by almost norming pairs in \cite[Section~5]{KMMPQ}.
\item If $G\in\Lin(X,Y)$ is a norm-one linear operator, then
\[
\nu_{L,G}(T)=\nu_G(T)
\quad
(T\in\Lin(X,Y)).
\]
Consequently, $n_{L,G}(X,Y)\leq n_G(X,Y).$
\end{enumerate}
\end{remark}

The connection with the linear theory is provided by Lipschitz-free spaces. For a scalar field $\mathbb{K}$, the \emph{Lipschitz-free space} $\F_\mathbb{K}(M)$ over a metric space $M$ is the closed
linear span in $\Lipz(M,\K)^*$ of the evaluation functionals
$\delta_p$, $p\in M$. For $p\neq q$, the elements
\[
m_{p,q}
=
\frac{\delta_p-\delta_q}{d(p,q)} \in S_{\F_{\mathbb{K}}(M)}
\]
are called \emph{molecules}. If $Y$ is a Banach space over $\mathbb{K}$, every $f\in\Lipz(M,Y)$ admits a unique
linearization $T_f\in\Lin(\F_\mathbb{K}(M),Y)$ such that
\[
T_f(\delta_p)=f(p)
\quad
(p\in M).
\]
The correspondence $f\mapsto T_f$ is a surjective $\mathbb K$-linear isometry from $\Lipz(M,Y)$ onto $\Lin(\F_{\mathbb K}(M),Y)$; see \cite{GK2003} for the real case, and \cite[Section~2]{ACP} and \cite[Section 5]{KMMPQ} for the complex case. In particular,
\[
\langle f,m_{p,q}\rangle :=  T_f(m_{p,q})  \in Y \quad\text{and}\quad \|T_f\|=\|f\|_L.
\]
Taking $Y=\mathbb K$ yields the canonical isometric identification
\[
\Lip_0 (M,\mathbb{K})=\mathcal{F}_\mathbb{K} (M)^*.
\]
We reserve the notation $\mathcal{F}(M)$ for the real Lipschitz-free space $\mathcal{F}_\mathbb{R} (M)$. 

\subsection{Main results}

Lipschitz-free linearization identifies the new numerical quantities
with their linear counterparts:
\[
\nu_{L,g}(f)=\nu_{T_g}(T_f)\quad\text{and}\quad n_{L,g}(M,Y)=n_{T_g}(\F_{\mathbb{K}}(M),Y).
\] Thus, the $g$-Lipschitz numerical radius and index may be viewed either as abstract numerical quantities on $\Lipz(M,Y)$ or as numerical quantities with respect to the linearization $T_g$. These identities, established in Proposition~\ref{prop:abstract-range}, provide the starting point for the following results.

Our first main result describes the real-valued theory in metric
terms. We state its finite version here; a more general formulation
appears in Theorem~\ref{thm:n-L-g-metric}.

\begin{maintheorem}\label{thm:intro-metric}
Let $M$ be a finite pointed metric space, and let $g\in S_{\Lipz(M,\R)}$. Let $\rho_g(p,q)$ be the infimum of the lengths of chains from $p$ to $q$ whose steps $(u,v)$ satisfy $|g(u)-g(v)|=d(u,v)$. Then
\[
\nu_{L,g}(f)
=\sup_{p\ne q}\frac{|f(p)-f(q)|}{\rho_g(p,q)}
\]
for $f\in\Lipz(M,\R)$ and
\[
n_{L,g}(M,\R)
=\min_{p\ne q}\frac{d(p,q)}{\rho_g(p,q)},
\]
where $\inf\emptyset=\infty$ and $\tfrac{a}{\infty}=0$ for finite $a$.
Moreover, $n_{L,g}(M,\R)>0$ if and only if $g$ is an extreme
point of $B_{\Lipz(M,\R)}$.
\end{maintheorem}

For shortest-path metrics on graphs, the metric description can be expressed in terms of the subgraph formed by the edges on which $g$ attains its Lipschitz norm. This yields an explicit formula for $n_{L,g}$ and the only possible positive values of the index are $1$ and the reciprocals of even positive integers. The same analysis also gives consequences for the classical numerical indices of the associated Lipschitz-free space and its dual.

\begin{maintheorem}\label{thm:intro-graphs}
Let $\Gamma=(V,E)$ be a connected simple unweighted graph with at least two vertices, equipped with its shortest-path metric and a distinguished vertex. 
\begin{enumerate}
\item[\textup{(a)}] For $g\in S_{\Lipz(V,\R)}$, let
\[
H_g=(V,\{uv\in E:|g(u)-g(v)|=1\}) \quad \text{and} \quad  D_g=\sup_{uv\in E}d_{H_g}(u,v),
\]
where vertices in different components of $H_g$ have distance
$\infty$. Then
\[
n_{L,g}(V,\R)=\frac{1}{D_g} \in\{0,1\}\cup\left\{\frac{1}{2k}:k\in\N\right\}.
\]
\item[\textup{(b)}] The classical numerical indices satisfy
\[
\frac12\le n(\Lipz(V,\R))\le n(\F(V))\le1.
\]
Moreover, $\Gamma$ is bipartite if and only if $n(\Lipz(V,\R)) =1$ if and only if  $n(\F(V))=1$.
\end{enumerate}
\end{maintheorem}

Theorem~\ref{thm:intro-graphs}(b) shows that either the Lipschitz-free space or its dual detects bipartiteness through its classical numerical index. The characterization applies to both finite and infinite connected graphs, and the universal lower bound $\tfrac12$ is sharp, as witnessed by the triangle. We also obtain exact formulas for cycle graphs. In particular, odd cycles exhibit markedly different behavior for the relative Lipschitz index and the classical numerical index:
\[
\sup_{g\in S_{\Lip_0 (C_{2k+1},\mathbb R)}}
n_{L,g}(C_{2k+1},\mathbb R) =\frac{1}{2k}\to 0, \quad n(\mathcal F(C_{2k+1})) =\frac{2k-1}{2k}\to 1.
\]

The metric structure of the domain also governs numerical-radius attainment. Earlier work on strongly norm-attaining Lipschitz functions revealed an obstruction to density arising from metric convexity and the length-property. Kadets, Mart\'in, and Soloviova \cite[Theorem 2.3]{KMS} proved non-density on metrically convex spaces, and Cascales et al. extended this phenomenon to complete length spaces \cite[Theorem 2.2]{CCGMR}. In the numerical-radius setting, Choi, Jung, and Tag \cite[Theorem 2.4]{CJT} showed that Lipschitz numerical-radius attaining self-maps are not dense on any Banach space. 

We extend the latter result to arbitrary norm-one reference maps and arbitrary Banach targets on length spaces. Let $\LipNRA_g (M,Y)$ denote the set of maps attaining their $g$-Lipschitz numerical radius (see Definition \ref{def:g-Lipschitz-numerical-radius-attaining}).

\begin{maintheorem}\label{thm:characterization-of-length}
Let $M$ be a complete pointed metric space and $Y$ be a real or complex Banach space. $M$ is a length space if and only if 
\[
\overline{\LipNRA_g (M,Y)} \neq \Lip_0 (M,Y) \quad \text{for every}~g \in S_{\Lip_0 (M,Y)}.
\]
Moreover, if $M$ is not a length space, one can in fact choose a single reference map $g$ for which every $f\in\Lip_0(M,Y)$ attains its $g$-Lipschitz numerical radius. 
\end{maintheorem}

Taking $M=Y=X$ and $g =\id_X$ recovers the non-density theorem of \cite{CJT}. We also characterize finite metric spaces: $M$ is finite if and only if every real-valued Lipschitz function attains its $g$-Lipschitz numerical radius for every norm-one Lipschitz map $g$.

Having considered the metric structure of the domain, we turn to the geometry of the Banach spaces involved. The next theorem gives upper bounds for maps taking values in $\ell_p(\Lambda)$ and lower bounds for the Lipschitz numerical index of spaces with property $\beta$.

\begin{maintheorem}\label{thm:intro-geometry}
\begin{enumerate}
\item[\textup{(a)}] Let $1<p<\infty$, and let $\Lambda$ be a set with at least two elements. For every pointed metric space $M$ and every $g\in S_{\Lipz(M,\ell_p(\Lambda))}$, where $\ell_p(\Lambda)$ is real or complex, we have
\[
n_{L,g}(M,\ell_p(\Lambda))
\le M_p
:=\max_{0\le t\le1}\frac{|t^{p-1}-t|}{1+t^p}.
\]
In particular, if $H$ is a real or complex Hilbert space
with $\dim H\ge2$, then
\[
n_{L,g}(M,H)=0
\]
for every pointed metric space $M$ and every $g\in S_{\Lipz(M,H)}$.
\item[\textup{(b)}] If a Banach space $X$ has property $\beta$ with constant $0\le\rho<1$, then
\[
n_L(X)\ge
\begin{cases}
\dfrac{1-\rho}{1+\rho}, & X\text{ is real},\\[1ex]
1-\rho, & X\text{ is complex}.
\end{cases}
\]
\end{enumerate}
\end{maintheorem}

For $H=\ell_2^2 (\mathbb{C})$ and $g = \text{Id}_H$, this recovers the example $n_L (H)=0<n(H)=\tfrac12$ discussed in \cite{PH}. A consequence of Theorem \ref{thm:intro-geometry}(a) establishes the same vanishing conclusion for every metric domain and every norm-one reference map into any real or complex Hilbert space of dimension at least two. 
The bounds in Theorem~\ref{thm:intro-geometry}(b) extend the corresponding classical estimates to the Lipschitz setting and answer \cite[Question~4.11]{CJT}. These bounds need not persist when the reference operator is perturbed: if $X$ has property $\beta$ and $\dim X\ge2$, Proposition~\ref{prop:property-beta-instability} provides norm-one operators $G$ arbitrarily close to $\id_X$ with $n_{L,G}(X,X)=n_G(X,X)=0$.

Finally, we compare the Lipschitz numerical index with the numerical indices of a linear reference operator and its adjoint. For arbitrary real Banach spaces, we obtain a lower bound in terms of the adjoint that complements the immediate upper bound $n_{L,G}(X,Y) \le n_G(X,Y)$. We then identify broad classes for which the upper bound is attained. The equality results extend those of P\'erez-Hern\'andez \cite{PH} from the identity to arbitrary norm-one operators and from dual targets to spaces that are $1$-complemented in their biduals.

\begin{maintheorem}\label{thm:intro-comparison}
Let $X$ and $Y$ be real Banach spaces, and let $G\in S_{\Lin(X,Y)}$.
\begin{enumerate}
    \itemsep0.25em 
    \item[\textup{(a)}] For every $f \in \Lip_0 (X,Y)$, 
    \[
    n_{G^*} (Y^*,X^*)\|f\|_L \le \nu_{L,G}(f).
    \]
    Consequently, 
    \[
    n_{G^*} (Y^*,X^*) \le n_{L,G} (X,Y) \le n_{G}(X,Y) .
    \]
    \item[\textup{(b)}] If either $X$ is separable or $Y$ is $1$-complemented in its bidual, then 
    \[
    n_{L,G}(X,Y)=n_G(X,Y).
    \]
\end{enumerate}
\end{maintheorem}

In particular, part (a) gives 
\[
n(X^*)\le n_L (X) \le n(X)
\]
for every real Banach space $X$, without separability or complementation assumptions. Hence $n_L(X)=n(X)$ whenever $n(X^*)=n(X)$.

Part (b) extends the equality results of P\'erez-Hern\'andez \cite[Theorems~2.1 and~3.1]{PH} and separates the roles of the domain and the target. Our operator formulation shows that separability of the domain alone suffices for every real Banach target, while $1$-complementation of the target in its bidual alone suffices for every real Banach domain. For example, $n_{L,G}(X,L_1[0,1])=n_G(X,L_1[0,1])$ for every real Banach space $X$ and every norm-one operator $G\colon X\to L_1[0,1]$. Moreover, taking $G=\id_X$ yields an extension even in the original self-map setting: the equality $n_L(X)=n(X)$ holds for every real Banach space that is $1$-complemented in its bidual, extending the previously known result for dual spaces.

\section{Abstract numerical ranges and first consequences}\label{sec:prelim}

\subsection{Abstract numerical ranges and linearization}

Given a Banach space $E$ and $u \in S_E$, the \emph{abstract numerical range of $z \in E$ with respect to $u$} is
\[
V(E,u,z) = \left\{ \varphi(z): \varphi\in S_{E^*},\ \varphi(u)=1 \right\}.
\]
The corresponding \emph{numerical radius} is defined by
\[
\nu(E,u,z) = \sup\left\{ |\lambda|:\lambda\in V(E,u,z) \right\}
= \sup\left\{ |\varphi(z)|: \varphi\in S_{E^*},\ \varphi(u)=1 \right\},
\]
and the \emph{numerical index of $E$ with respect to $u$} is 
\[
n(E, u)= \inf \{ \nu (E,u,z): z \in S_E \}. 
\]

Recall that if $G$ is a norm-one operator between Banach spaces $X$ and $Y$, then the numerical radius $\nu_G(T)$ coincides with the abstract numerical radius $\nu(\Lin(X,Y),G,T)$ for every $T\in\Lin(X,Y)$; see \cite[Section~1]{KMMPQ}. Consequently, $n_G(X,Y)=n(\Lin(X,Y),G)$. We shall use these identities without further mention.

\begin{proposition}\label{prop:abstract-range}
    Let $M$ be a pointed metric space, $Y$ be a Banach space, and let $f,g\in\Lipz(M,Y)$ satisfy $\|g\|_L=1$. Then 
    \begin{enumerate}
        \item[\textup{(a)}] $\overline{\conv} W_{L,g} (f) = V(\Lip_0 (M,Y), g,f) = V(\mathcal{L}(\mathcal{F}_{\mathbb{K}} (M), Y), T_g, T_f)$,
        \item[\textup{(b)}] $\nu_{L,g}(f)=\nu(\Lip_0(M,Y),g,f)=\nu_{T_g}(T_f) = \max \{|\lambda|:\lambda \in W_{L,g} (f) \}$, 
        \item[\textup{(c)}] $n_{L,g} (M,Y)=n(\Lip_0 (M,Y), g)=n_{T_g}(\mathcal{F}_{\mathbb{K}} (M),Y)$.
    \end{enumerate}
\end{proposition}

\begin{proof}
Consider the following $\mathbb{T}$-invariant subset of $B_{\Lip_0(M,Y)^*}$:
\[
\mathcal A
=
\left\{
\Phi_{p,q,y^*}: p,q\in M,\ p\neq q,\ y^*\in S_{Y^*} \right\},
\]
where
\[
\Phi_{p,q,y^*}(h) = y^* \left( \frac{h(p)-h(q)}{d(p,q)} \right)
\]
for $h\in \Lipz(M,Y)$. Clearly, $\mathcal A\subseteq B_{\Lip_0(M,Y)^*}$, and for every $h\in \Lip_0(M,Y)$,
\[
\|h\|_L = \sup_{p\neq q} \left\|\frac{h(p)-h(q)}{d(p,q)} \right\| = \sup_{\Phi\in\mathcal A}|\Phi(h)| = \sup_{\Phi\in\mathcal{A}} \Rea\Phi(h).
\]
Consequently, 
$B_{\Lip_0(M,Y)^*}=\overline{\conv}^{\,w^*}(\mathcal A)$. By \cite[Proposition 2.14]{KMMPQ}, 
\begin{align*}
V(\Lip_0(M,Y),g,f) = \overline{\conv}W_{L,g}(f), \quad \text{and} \quad \nu(\Lip_0(M,Y),g,f)=\nu_{L,g}(f).
\end{align*}
The sets $\overline{W_{L,g}^\delta(f)}$ are nonempty compact sets and decrease as $\delta \downarrow0$. Thus, 
\[
\nu_{L,g}(f)=\max \{|\lambda|:\lambda \in W_{L,g}(f)\}.
\]
Finally, the remaining identities follow from the surjective linear isometry 
\[
\Lip_0(M,Y)\to \mathcal L(\mathcal{F}_{\mathbb{K}} (M),Y), \quad h \mapsto T_h. \qedhere
\]
\end{proof}

\subsection{Norm formulas and basic examples}

Recall that a point $u\in S_E$ of a Banach space $E$ is called a \emph{strongly extreme point} of $B_E$ if, for any sequences $(u_n)$ and $(v_n)$ in $B_E$,
\[
    \left\|u-\frac{u_n+v_n}{2}\right\|\to 0 \implies  \|u_n-v_n\|\to 0.
\]
Every strongly extreme point of $B_E$ is an extreme point of $B_E$.

\begin{corollary}\label{cor:norm-derivative}
    Let $M$ be a pointed metric space, $Y$ be a Banach space, and let $g\in S_{\Lipz(M,Y)}$. Then the following assertions hold: 
    \begin{enumerate}
        \item[\textup{(a)}] For every $f\in\Lipz(M,Y)$,
\begin{align*}
\vL{g}(f)=\max_{\theta\in\T} \lim_{\alpha\to0^+}\frac{\|g+\alpha\theta f\|_L-1}{\alpha} = \lim_{\alpha \to 0^+} \max_{\theta\in\mathbb{T}} \frac{\|g+\alpha\theta f\|_L-1}{\alpha}. 
\end{align*}
\item[\textup{(b)}]     For $0<\lambda \leq 1$, the inequality $n_{L,g}(M,Y)\geq \lambda$ holds if and only if 
             \[
             \max_{\theta \in \mathbb{T}} \| g+\theta f \|_L \geq 1 + \lambda \|f\|_L \quad \text{for every}~f \in \Lip_0(M,Y).
             \]
    In particular, $n_{L,g}(M,Y)=1$ if and only if $g$ is a spear element of $\Lip_0(M,Y)$. 
    \item[\textup{(c)}] If $n_{L,g}(M,Y)>0$, then $g$ is a strongly extreme point of $B_{\Lip_0 (M,Y)}$.    
    \end{enumerate}
\end{corollary}

\begin{proof} Apply Proposition \ref{prop:abstract-range} together with \cite[Lemmas 2.2 and 2.3, and Proposition 2.5]{KMMPQ} to $E=\Lip_0(M,Y)$. 
\end{proof}

We next use the linearization identity in Proposition~\ref{prop:abstract-range} to compute the index for rank-one Lipschitz maps.

\begin{corollary}\label{cor:rank-one}
    Let $M$ be a pointed metric space, $Y$ be a Banach space, and let $\phi \in \Lip_0 (M,\mathbb{K})$ and $y_0 \in Y$ such that $\|\phi\|_L=\|y_0\| =1$. If $g= \phi \otimes y_0$, then 
    \[
    n_{L,g} (M,Y) = n(\Lip_0 (M,\mathbb{K}), \phi) n(Y,y_0).
    \]
\end{corollary}

\begin{proof}
By Proposition~\ref{prop:abstract-range},
\[
n_{L,g}(M,Y) = n_{T_g}\bigl(\mathcal F_{\mathbb{K}}(M),Y\bigr).
\]
Since $g=\phi\otimes y_0$, its linearization is the rank-one operator $T_g=T_\phi\otimes y_0$, where $T_\phi\in\mathcal F_{\mathbb{K}}(M)^*$ corresponds canonically to $\phi\in\Lip_0(M,\mathbb{K})$. Hence, by \cite[Proposition~3.6]{KMMPQ},
\begin{align*}
n_{L,g}(M,Y) &=
n_{T_\phi\otimes y_0}\bigl(\mathcal F_{\mathbb{K}}(M),Y\bigr) \\
&=
n\bigl(\mathcal F_{\mathbb{K}}(M)^*,T_\phi\bigr)n(Y,y_0)= n\bigl(\Lip_0(M,\mathbb{K}),\phi\bigr)n(Y,y_0).\qedhere
\end{align*}
\end{proof}

\begin{corollary}\label{cor:interval}
Let $I \subseteq \mathbb{R}$ be an interval containing $0$ and $g \in S_{\Lip_0 (I,\mathbb{K})}$. Then 
\[
n_{L,g} (I,\mathbb{K})=
\begin{cases}
    1,  \quad |g'(t)|=1 &\text{for a.e. $t \in I$},\\ 
    0, \quad &\text{otherwise}.
\end{cases}
\]
\end{corollary}

\begin{proof}
    By Corollary \ref{cor:rank-one}, we have $n_{L,g}(I,\mathbb{K})= n(\Lip_0 (I,\mathbb{K}), g)$. Note that $\Phi: \Lip_0 (I,\mathbb{K}) \to L_\infty(I,\mathbb{K})$, $\Phi(f)=f'$, is a surjective linear isometry. It follows from \cite[Lemma 2.4]{KMMPQ} that 
    \[
    n(\Lip_0 (I,\mathbb{K}),g) = n (\Phi(\Lip_0 (I,\mathbb{K})), \Phi(g)) =  n (L_\infty (I,\mathbb{K}), g'). 
    \]
    If $|g'(t)|=1$ for a.e. $t \in I$, then $g'$ is a spear element of $L_\infty (I,\mathbb{K})$ (\cite[Example 2.12(d)]{KMMP}), so $n(L_\infty (I,\mathbb{K}),g')=1$. Otherwise, $g'$ is not an extreme point of $B_{L_\infty (I,\mathbb{K})}$; hence \cite[Lemma 2.3]{KMMPQ} shows that $n(L_\infty (I,\mathbb{K}), g')=0$. 
\end{proof}

\subsection{Weak G\^ateaux derivative formulas}

We conclude this section with a description of the $g$-Lipschitz numerical radius in terms of weak G\^ateaux derivatives when the domain is a Banach space. Recall that a map $f:X\to Y$ between Banach spaces is weakly G\^ateaux differentiable at $x_0\in X$ if there is $df(x_0)\in\Lin(X,Y)$ such that
\[
y^*\left(\frac{f(x_0+tu)-f(x_0)}{t}\right)\to y^*(df(x_0)u)\quad\text{as}\quad t\to0~\text{in}~\R,
\]
for all $u\in X$ and $y^*\in Y^*$. In this case, the operator $df(x_0)$ is unique and is called the weak G\^ateaux derivative of $f$ at $x_0$.

Part (a) of Proposition~\ref{prop:weak-derivative-formulas} will be used in Section~\ref{sec:comparison} to compare the Lipschitz and linear numerical indices associated with a linear reference operator.

\begin{proposition}\label{prop:weak-derivative-formulas}
    Let $X$ and $Y$ be Banach spaces, and let $f,g\in S_{\Lip_0(X,Y)}$. 
    \begin{enumerate}
        \item[\textup{(a)}] Suppose that $f,g \in \Lip_0 (X,Y)$ are weakly G\^ateaux differentiable at $x_0 \in X$ and $\|dg(x_0)\|=1$. Then 
        \[
            \nu_{dg(x_0)} (df(x_0) ) \leq \nu_{L,g} (f).
        \]
        \item[\textup{(b)}] Suppose that $f,g \in \Lip_0 (X,Y)$ are weakly G\^ateaux differentiable at every point of $X$. For $0<\delta<1$, set 
                \begin{align*}\label{def-v-delta-f-g}
            &\mathfrak{v}_\delta(f,g) := \sup \left \{|y^* ( df(x) (u) ) | : \begin{aligned} &x \in X, \, u\in S_X, \, y^* \in S_{Y^*}, \\
            &\Rea y^*(dg(x)(u) )>1-\delta \end{aligned} \right \}.
        \end{align*}  Then 
        \begin{equation*}
            \nu_{L,g} (f) = \inf_{0<\delta<1} \mathfrak{v}_\delta(f,g).
         \end{equation*} 
         In particular, if $g=G\in S_{\mathcal{L}(X,Y)}$, then 
         \[
             \nu_{L,G} (f) = \inf_{0<\delta<1} \sup \{ \nu_{G}^\delta  (df(x)) : x \in X \}. 
         \]
    \end{enumerate}
\end{proposition}

\begin{proof}
    (a) Fix $0<\delta<1$, $u \in S_X$, and $y^*\in S_{Y^*}$ such that 
    \[
    \Rea y^*(dg(x_0)u)>1-\delta. 
    \]
    For all sufficiently small $t>0$, we have 
    \[
    \Rea y^* \left( \frac{g(x_0+tu)-g(x_0)}{t}\right) >1-\delta.
    \]
    Hence 
    \[
    \left| y^* \left( \frac{f(x_0+tu)-f(x_0)}{t} \right) \right| \leq \nu_{L,g}^\delta (f).
    \]
    Letting $t\downarrow0$, we have 
    \[
    |y^* (df(x_0)(u) )| \leq \nu_{L,g}^\delta (f).
    \]
    Consequently, $\nu_{dg(x_0)}^\delta (df(x_0))\leq \nu_{L,g}^\delta (f)$. Letting $\delta \to 0$, we conclude that $\nu_{dg(x_0)} (df(x_0) ) \leq \nu_{L,g} (f)$.

    (b) First observe that $\|h\|_L = \sup_{x \in X} \|dh(x)\|$ for every $h \in \Lip_0 (X,Y)$ that is weakly G\^ateaux differentiable at every point of $X$. In particular, we have 
    \[
    \|h\|_L = \sup_{x \in X} \|dh(x)\| \quad \text{for every}~h \in E:= \spann_{\mathbb{K}} \{ f,g\}.
    \]
    Define $\Psi_{x,u,y^*} \in E^*$ by $\Psi_{x,u,y^*} (h)=y^*(dh(x)u)$, and $\mathcal{A}=\{\Psi_{x,u,y^*}: x\in X, \, u\in S_X,\, y^* \in S_{Y^*}\}$. Then $\mathcal{A}\subseteq B_{E^*}$, and 
    \[
    \|h\|_L = \sup_{\Psi \in \mathcal{A}} \Rea\Psi(h)
    \]
    for $h\in E$. Thus, $B_{E^*}=\overline{\text{conv}}^{w^*}(\mathcal{A})$, and \cite[Proposition 2.14]{KMMPQ} implies that 
    \[
    \nu(E,g,f)=\inf_{0<\delta<1} \mathfrak{v}_\delta(f,g).
    \]
    Finally, since the inclusion $\iota\colon E\to\Lip_0(X,Y)$ is a linear isometric embedding, applying \cite[Lemma 2.4(b)]{KMMPQ} to $\iota$ together with Proposition \ref{prop:abstract-range}, we obtain 
    \[
    \nu(E,g,f)=\nu(\Lip_0 (X,Y),g,f)=\nu_{L,g}(f).
    \]
    This completes the proof.
\end{proof}

\section{Lipschitz numerical radius: metric interpretation}\label{sec:metric-interpretation}

Throughout this section, we consider \emph{real-valued} Lipschitz functions, \emph{real} Lipschitz-free spaces, and adopt the conventions
\[
\inf\emptyset=\infty, \quad\text{and}\quad \frac{a}{\infty}=0
\]
for every $a\in\R$.

\subsection{Almost-norming path metrics}
For $0<\delta<1$ and $g \in S_{\Lip_0 (M,\mathbb{R})}$, consider 
\[
\mathfrak{N}_\delta (g) := \left \{ (p,q)\in M\times M: p\neq q, \, \frac{|g(p)-g(q)|}{d(p,q)}>1-\delta \right \}.
\]
For $f \in \Lip_0 (M,\mathbb{R})$, notice that 
\begin{equation}\label{eq:nu-L-g-f-mathfrak-N-description}
\nu_{L,g}^\delta (f) = \sup \left \{ \frac{|f(p)-f(q)|}{d(p,q)} : p\neq q,\, (p,q) \in \mathfrak{N}_\delta (g) \right \} .    
\end{equation}

Let us consider the following extended metric $\rho_{g,\delta}$ on $M$:
\begin{equation}\label{def-rho-g-delta}
\rho_{g,\delta} (p,q) := 
\inf\left\{
\sum_{j=1}^n d(x_j,x_{j-1}):
\begin{aligned}
&p=x_0,x_1,\ldots,x_n=q,\\
&(x_j,x_{j-1})\in\mathfrak N_\delta (g), ~ j=1,\ldots, n
\end{aligned}
\right\}.
\end{equation}
If there is no such chain $\{x_0,\ldots,x_n\}$ in \eqref{def-rho-g-delta} between $p$ and $q$, then we put $\rho_{g,\delta}(p,q)=+\infty$; for $p=q$ we allow the trivial chain, so that $\rho_{g,\delta}(p,p)=0$. Lemma~\ref{lem:basic-property-of-rho_g_delta} is easy to verify, so we omit the proof.

\begin{lemma}\label{lem:basic-property-of-rho_g_delta}
Let $g \in S_{\Lip_0 (M,\mathbb{R})}$ and $0<\delta<1$ be given. Then 
\begin{enumerate}
    \item[\textup{(a)}] $d(p,q)\leq \rho_{g,\delta} (p,q)$ for all $p,q\in M$;
    \item[\textup{(b)}] if $(p,q)\in\mathfrak{N}_\delta(g)$, then $d(p,q)=\rho_{g,\delta}(p,q)$;
    \item[\textup{(c)}] if $0<\delta_1\le \delta_2<1$, then $\rho_{g,\delta_1}(p,q) \ge \rho_{g,\delta_2}(p,q)$.
\end{enumerate}
\end{lemma}

\begin{proposition}\label{prop:nu-L-g-delta-description} Let $M$ be a metric space, $g \in S_{\Lip_0 (M,\mathbb{R})}$, $f \in \Lip_0 (M,\mathbb{R})$ and $0<\delta<1$ be given. Then 
    \[ 
    \nu_{L,g}^\delta (f) = \sup \left\{ \frac{|f(p)-f(q)|}{\rho_{g,\delta}(p,q)} : p\neq q \in M\right\}.
    \]
\end{proposition}

\begin{proof}
    For simplicity, let us denote by $A$ the right-hand side term in \eqref{eq:nu-L-g-f-mathfrak-N-description}. First, if $p\neq q\in M$ and $(p,q) \in \mathfrak{N}_\delta(g)$, then $d(p,q)=\rho_{g,\delta}(p,q)$. Thus, 
    \[
    \frac{|f(p)-f(q)|}{d(p,q)} = \frac{|f(p)-f(q)|}{\rho_{g,\delta}(p,q)} \leq \sup \left\{ \frac{|f(x)-f(y)|}{\rho_{g,\delta}(x,y)} : x\neq y \in M\right\} =:B.
    \]
    This proves that $A\leq B$. 

    For the reversed inequality, let $x \neq y$ be given. We claim that 
    \[
    \frac{|f(x)-f(y)|}{\rho_{g,\delta}(x,y)} \leq A.
    \]
    If $\rho_{g,\delta} (x,y) =+\infty$, then there is nothing to prove. Let $\{x=x_0,x_1,\ldots, x_n=y\}$ be a chain from $x$ to $y$ such that $(x_j,x_{j-1}) \in \mathfrak{N}_\delta (g)$. Then 
    \begin{align*}
        |f(x)-f(y)| &\leq \sum_{j=1}^n |f(x_j)-f(x_{j-1})| \leq A\sum_{j=1}^n d(x_j,x_{j-1}) .
    \end{align*}
    This implies that $|f(x)-f(y)| \leq A \rho_{g,\delta}(x,y)$ and the claim is proved. Taking the supremum over $x\neq y \in M$, we obtain that $B \leq A$. 
\end{proof}

Writing $\|f\|_{L,\rho_{g,\delta}}$ for the Lipschitz constant
of $f$ with respect to the extended metric $\rho_{g,\delta}$,
Proposition~\ref{prop:nu-L-g-delta-description} gives
\[
\nu_{L,g}^{\delta}(f)
=
\|f\|_{L,\rho_{g,\delta}}.
\]

Given $g \in S_{\Lip_0 (M,\mathbb{R})}$, consider the following set and value as well:
\begin{align*}
&\mathfrak{N} (g) := \left \{ (p,q)\in M\times M: p\neq q, \, \frac{|g(p)-g(q)|}{d(p,q)}=1 \right \}, \quad\text{and}\\
&\rho_g(p,q):=\inf\left\{
\sum_{j=1}^n d(x_j,x_{j-1}):
\begin{aligned}
&p=x_0,x_1,\ldots,x_n=q,\\
&(x_j,x_{j-1})\in\mathfrak N(g),~ j=1,\ldots, n
\end{aligned}
\right\}.\end{align*}
The set $\mathfrak{N} (g)$ is possibly empty, and in this case, $\rho_g(p,q) = \infty$ for $p\neq q$. Nevertheless, if $M$ is a finite metric space, then we observe the following.

\begin{corollary}\label{cor:nu-L-g-description-finite-metric}
    Let $M$ be a finite metric space, $g \in S_{\Lip_0 (M,\mathbb{R})}$, and $f \in \Lip_0 (M,\mathbb{R})$. Then 
     \[ 
    \nu_{L,g} (f) = \sup \left\{ \frac{|f(p)-f(q)|}{\rho_{g}(p,q)} : p\neq q \in M\right\}.
    \]
\end{corollary}

\begin{proof} Observe that for sufficiently small $\delta$, we have $\mathfrak{N}(g)=\mathfrak{N}_\delta (g)$ and $\rho_{g,\delta}=\rho_g$. Thus, the result follows from Proposition \ref{prop:nu-L-g-delta-description}.
\end{proof}

\subsection{Metric distortion and the numerical index}

For $0<\delta<1$, define
\[
M_g^\delta:= \sup \left\{ \frac{\rho_{g,\delta} (p,q) }{d(p,q)} : p\neq q \in M \right\} \in [1,\infty].
\]

\begin{theorem}\label{thm:n-L-g-metric} 
    Let $M$ be a metric space and let $g \in S_{\Lip_0 (M,\mathbb{R})}$. For every $0<\delta<1$ such that either $M$ is uniformly discrete or $M_g^\delta <\infty$, we have 
    \[
    \inf \{ \nu_{L,g}^\delta (f):f\in\Lip_0(M,\mathbb{R}), \,\|f\|_L=1\} = \frac{1}{M_{g}^\delta}.
    \]
    Consequently, if either $M$ is uniformly discrete or $M_g^\delta <\infty$ for every $0<\delta<1$, then 
    \begin{equation}\label{thm:n-L-g-metric-1}
    n_{L,g} (M,\mathbb{R}) = \inf_{0<\delta<1} \frac{1}{M_g^\delta}.
    \end{equation}
    In particular, if $M$ is finite, then 
    \begin{equation}\label{thm:n-L-g-metric-2}
    n_{L,g}(M,\mathbb{R})=\min \left\{ \frac{d(p,q)}{\rho_g(p,q)} : p\neq q \in M \right\}.
    \end{equation}
\end{theorem}

\begin{proof}
Fix $0<\delta<1$ and write 
\[
a_\delta =  \inf \{ \nu_{L,g}^\delta (f):f\in\Lip_0(M,\mathbb{R}), \,\|f\|_L=1\}.
\]
Suppose that $M_g^\delta<\infty$. By Proposition \ref{prop:nu-L-g-delta-description}, we have for $f \in \Lip_0(M,\mathbb{R})$ and $p\neq q \in M$ that 
\[
|f(p)-f(q)| \leq \nu_{L,g}^\delta (f) \rho_{g,\delta} (p,q) \leq \nu_{L,g}^\delta(f) M_g^\delta d(p,q).
\]
Thus, $\|f\|_L \leq \nu_{L,g}^\delta(f) M_{g}^\delta$, which implies that $a_\delta \geq \tfrac{1}{M_g^\delta}$. 

For the reverse inequality, fix $p\neq q$ and define 
\[
h_q (u):= \rho_{g,\delta} (u,q)-\rho_{g,\delta} (0,q)
\]
for every $u\in M$. Note that for $u\neq v \in M$, 
\begin{equation}\label{eq:hqu-hqv}
|h_q(u)-h_q(v)| = |\rho_{g,\delta} (u,q)-\rho_{g,\delta}(v,q)|\leq \rho_{g,\delta}(u,v)\leq M_g^\delta d(u,v).
\end{equation}
This shows that $h_q \in \Lip_0 (M,\mathbb{R})$. By Proposition \ref{prop:nu-L-g-delta-description} and \eqref{eq:hqu-hqv}, we have $\nu_{L,g}^\delta (h_q)\leq 1$. On the other hand, 
\[
\|h_q\|_L \geq \frac{|h_q(p)-h_q(q)|}{d(p,q)} = \frac{\rho_{g,\delta}(p,q)}{d(p,q)}>0.
\]
It follows that 
\[
a_\delta \leq \nu_{L,g}^\delta \left( \frac{h_q}{\|h_q\|_L} \right) \leq \frac{1}{\|h_q\|_L} \leq \frac{d(p,q)}{\rho_{g,\delta}(p,q)}.
\]
Since $p\neq q \in M$ was arbitrary, we conclude that 
\[
a_\delta \leq \inf \left\{ \frac{d(p,q)}{\rho_{g,\delta}(p,q)}:p\neq q \in M\right\} = \frac{1}{M_g^\delta}.
\]

Finally, suppose that $M$ is uniformly discrete and $M_g^\delta=\infty$. Let $\sigma := \inf \{d(p,q):p\neq q\in M\} >0$. For $q \in M$ and $R>0$, define $h_{q,R}:M\to\mathbb{R}$ by  
\[
h_{q,R} (u) = \min \{ \rho_{g,\delta} (u,q),R\} - \min \{\rho_{g,\delta } (0,q),R\}
\]
for every $u\in M$. Here, we understand $\min \{\infty, R\}=R$. Note that $\|h_{q,R}\|_L\leq \tfrac{R}{\sigma}$. Moreover, for $u \neq v \in M$, we have that 
\begin{itemize}
    \item if $\rho_{g,\delta} (u,q), \rho_{g,\delta}(v,q)<\infty$, then 
    \[
    |h_{q,R}(u)-h_{q,R}(v)| \leq |\rho_{g,\delta}(u,q)-\rho_{g,\delta}(v,q)| \leq \rho_{g,\delta} (u,v).
    \]
    \item if $\rho_{g,\delta} (u,q) =\infty$ and $\rho_{g,\delta}(v,q)<\infty$, then $\rho_{g,\delta}(u,v)=\infty$. Thus, 
    \[
    |h_{q,R}(u)-h_{q,R}(v)| < \infty = \rho_{g,\delta} (u,v).
    \]
    \item if $\rho_{g,\delta} (u,q)=\rho_{g,\delta}(v,q) =\infty$, then 
    \[
    |h_{q,R}(u)-h_{q,R}(v)| = |R-R|=0 \leq \rho_{g,\delta} (u,v).
    \]
\end{itemize}
Therefore, it follows from Proposition \ref{prop:nu-L-g-delta-description} that $\nu_{L,g}^\delta (h_{q,R})\leq 1.$
On the other hand,
\[
\|h_{q,R}\|_L \geq \left|\frac{h_{q,R}(p)-h_{q,R}(q)}{d(p,q)}\right| = \frac{\min \{ \rho_{g,\delta} (p,q), R\}}{d(p,q)},
\]
which implies that 
\[
a_\delta \leq \nu_{L,g}^\delta\left(  \frac{h_{q,R}}{\|h_{q,R}\|_L}\right) \leq \frac{1}{\|h_{q,R}\|_L} \leq \frac{d(p,q)}{\min \{\rho_{g,\delta}(p,q),R \}}. 
\]
Letting $R \to \infty$ and taking the infimum for $p\neq q \in M$, we obtain 
\[
a_\delta \leq \frac{1}{M_{g}^\delta }=0. \qedhere
\]
\end{proof}

The following example shows that formula \eqref{thm:n-L-g-metric-1} can fail without the hypotheses of Theorem \ref{thm:n-L-g-metric}.

\begin{example}\label{thm:n-L-g-metric-counterexample}
    Let $M=[0,1]$ and define $g \in \Lip_0(M,\mathbb{R})$ as 
    \[
    g(t)=\min\left\{t- \frac{1}{2^{n+1}}, \frac{1}{2^n}-t\right\} 
    \]
    for $t\in\left [\frac{1}{2^{n+1}},\frac{1}{2^n} \right]$, where $n \in \mathbb{N}\cup \{0\}$.
    Notice that $|g'(t)|=1$ for a.e. $t \in [0,1]$. By Corollary \ref{cor:interval}, we thus obtain that $n_{L,g}([0,1],\mathbb{R})=1.$ On the other hand, since $0\leq g(t)\leq \tfrac{t}{3}$ for all $t \in [0,1]$, we have that if $0<\delta<\tfrac23$, then $(0,t)\not\in \mathfrak{N}_\delta (g)$ for every $t >0$. It follows that $\rho_{g,\delta}(0,t)=\infty$; hence $M_g^\delta =\infty$. Therefore, \eqref{thm:n-L-g-metric-1} does not hold in this case.
\end{example}

\subsection{Finite metric spaces}

The following example shows that every value in $(0,1]$ can be achieved as the $g$-Lipschitz numerical index of a three-point metric space. 

\begin{example}\label{ex:three-point-example}
    Let $0<\tau \leq 1$. Consider $M_\tau =\{0,a,b\}$ with 
    \[
    d(a,b)=2\tau, \quad\text{and}\quad d(a,0)=d(b,0)=1.
    \]
    Consider $g: M_\tau \to \mathbb{R}$ given as 
    \[
    g(0)=0, \quad\text{and}\quad g(a)=g(b)=1.
    \]
    Then $\rho_g(a,b)=2$, and $\rho_g(a,0)=\rho_g(b,0)=1$. By Theorem \ref{thm:n-L-g-metric}, we have 
    \[
    n_{L,g}(M_\tau ,\mathbb{R})=\min\left\{ \frac{d(p,q)}{\rho_g (p,q)} : p \neq q \in M_\tau \right\} =\tau.
    \]

    On the other hand, the classical numerical index of $\mathcal{F}(M_\tau)$ can be computed from the formula \cite[Theorem 1.1]{CGM2024} as 
    \[ 
    n(\mathcal{F}(M_\tau))=\begin{cases} 
        \dfrac{1}{1+2\tau}, \quad &0<\tau\leq \dfrac{1}{2},\\
        \dfrac{1}{3-2\tau}, \quad &\dfrac{1}{2}\leq\tau\leq 1.
    \end{cases}
    \]
     For certain values of $\tau$, we have the following table. 
   \begin{table}[ht]
\centering
\renewcommand{\arraystretch}{1.3}
\begin{tabular}{*{3}{>{\centering\arraybackslash}p{2.6cm}}}
\toprule
$\tau$ 
& $n_{L,g}(M_\tau,\mathbb{R})$ 
& $n(\mathcal{F}(M_\tau))$ \\
\midrule
$1/4$ & $1/4$ & $2/3$ \\
$1/2$ & $1/2$ & $1/2$ \\
$3/4$ & $3/4$ & $2/3$ \\
\bottomrule
\end{tabular}
\end{table}

This shows that there is no inequality between $n_{L,g}(M,\mathbb{R})$ and $n(\mathcal{F}(M))$ in general. Also, notice that 
\[
n_{L,g}(M_{\tau},\mathbb{R})\to 0 \quad \text{while} \quad n(\mathcal{F}(M_\tau))\to 1
\]
as $\tau \downarrow0$.
\end{example}

The extreme points of the closed unit ball of $\Lip_0(M,Y)$ were characterized in \cite[Theorem~2.2]{BCJE} when $M$ is a finite pointed metric space and $Y$ is a strictly convex real Banach space. Specializing to $Y=\mathbb R$, their characterization can be expressed in our notation as follows: for $g\in S_{\Lip_0(M,\mathbb R)}$,
\[
    g\in \mathrm{ext}(B_{\Lip_0(M,\mathbb R)}) \iff    \rho_g(0,p)<\infty \,\,\text{ for every }p\in M.
\]
Together with Theorem~\ref{thm:n-L-g-metric}, this yields the
following characterization in terms of the $g$-Lipschitz
numerical index.
\begin{corollary}\label{cor:characterization-of-extreme-Lipschitz-functions-on-finite-metric-spaces}
    Let $M$ be a finite pointed metric space and let $g\in S_{\Lip_0(M,\mathbb R)}$. Then the following are equivalent. 
    \begin{enumerate}
        \item[\textup{(a)}] $n_{L,g}(M,\R)>0$;
        \item[\textup{(b)}] $\rho_{g}(p,q)<\infty$ for every $p,q \in M$; 
        \item[\textup{(c)}] $g \in \mathrm{ext}(B_{\Lipz(M,\mathbb{R})})$.
    \end{enumerate}
    Moreover, $n_{L,g}(M,\mathbb{R})=1$ if and only if $d=\rho_{g}$. 
\end{corollary}

Together with Corollary~\ref{cor:nu-L-g-description-finite-metric} and Theorem~\ref{thm:n-L-g-metric}, this characterization completes the proof of Theorem~\ref{thm:intro-metric}.

\section{Lipschitz numerical index of a graph}\label{sec:graph}

Throughout this section, a \emph{graph} is assumed to be a connected simple unweighted graph $\Gamma = (V,E)$ with at least two vertices. The set $V$ consists of its vertices and an element of $E$ is an unordered pair $\{u,v\}$ of distinct vertices. We also write $uv \in E$. We fix a distinguished vertex $0\in V$ and equip $V$ with the shortest-path metric $d=d_\Gamma$. All Lipschitz and Lipschitz-free spaces considered in this section are \emph{real} as in Section \ref{sec:metric-interpretation}.

\subsection{Saturated subgraphs and the index formula}

\begin{definition}
Let $\Gamma=(V,E)$ be a graph and let $g\in S_{\Lip_0 (V,\mathbb{R})}$. The set $E_g$ of \emph{$g$-saturated edges} is 
    \[
    E_g = \{uv\in E : |g(v)-g(u)|=d_\Gamma (u,v)\}
    \]
    and $ H_g = (V,E_g)$ is called a  \emph{$g$-saturated subgraph} of $\Gamma$.
\end{definition}

On every edge in $E_g$, the function $g$ attains its Lipschitz constant $1$.
The graph $H_g$ is obtained by deleting the unsaturated edges, without removing any vertices. We write $d_{H_g}$ for the shortest-path distance within $H_g$. If two vertices cannot be joined, then their distance is $\infty$. If $E_g=\emptyset$, then $H_g=(V,\emptyset)$ is a graph without edges, and $d_{H_g} (u,v) =\infty$ for every $u \neq v \in V$. Notice that we always have $d_\Gamma \le d_{H_g}$. Finally, we consider 
\[
D_g := \sup_{uv \in E } d_{H_g} (u,v).
\]

\begin{definition}
    Let $\Gamma=(V,E)$ be a graph and let $g\in S_{\Lip_0 (V,\mathbb{R})}$. For $0<\delta<1$, put 
    \[
    E^\delta_g = \{uv\in E: |g(v)-g(u)|>1-\delta\} 
    \]
    and write 
    \[
    a^\delta_g (f) = \sup \{|f(v)-f(u)|:uv \in E^\delta_g\}. 
    \]
\end{definition}

\begin{lemma}\label{lem:basic-lemma-on-graphs}
    Let $\Gamma=(V,E)$ be a graph, $g\in S_{\Lip_0 (V,\mathbb{R})}$, and let $f \in \Lip_0(V,\mathbb{R})$. Then
    \begin{enumerate}
        \item[\textup{(a)}] $\displaystyle \|f\|_L=\sup_{uv\in E} |f(v)-f(u)|$;
    \item[\textup{(b)}] for $0<\eta,\delta <1$, we have 
    \[
    a^\delta_g (f) \leq \nu_{L,g}^\delta (f) \leq a^\eta_g (f) + \frac{\delta}{\eta}\|f\|_L.
    \]
    Consequently, we have 
    \[
    \nu_{L,g}(f) = \lim_{\delta\to0^+} a^\delta_g (f).
    \]
    \item[\textup{(c)}] If $g$ is an extreme point of $B_{\Lip_0(V,\mathbb{R})}$, then $g(V)\subseteq \mathbb{Z}$.
        \end{enumerate} 
\end{lemma}

\begin{proof}
(a) It is enough to check that $\|f\|_L\leq \sup_{uv\in E} |f(v)-f(u)|.$ To this end, pick $p,q \in V$ and take $p=v_0,v_1,\ldots,v_r=q$ to be the shortest path connecting $p$ and $q$. Then 
\begin{align*}
|f(p)-f(q)| &\leq \sum_{j=1}^r |f(v_j)-f(v_{j-1})| \\
&\leq r \sup_{uv\in E}|f(v)-f(u)| = d(p,q) \sup_{uv\in E}|f(v)-f(u)|. 
\end{align*}
This completes the verification.

(b) It is clear by definition that $a^\delta_g(f) \leq \nu_{L,g}^\delta (f)$. To obtain the second inequality, let $p\neq q \in V$ be such that $|g(q)-g(p)|>(1-\delta)d(p,q)$. Let $p=v_0,v_1,\ldots,v_r=q$ be the shortest path connecting $p$ and $q$. Consider the sets
\begin{align*}
&A:= \{j \in \{1, \ldots, r\} : v_{j-1}v_j \in (E^\eta_g)^c\}, \\ 
&B:= \{j \in \{1, \ldots, r\} : v_{j-1}v_j \in E^\eta_g\},
\end{align*}
and suppose that $A$ has cardinality $s$. Observe that 
\begin{align*}
    (1-\delta)r < |g(q)-g(p)| &\leq \sum_{j=1}^r |g(v_j)-g(v_{j-1})| \\
    &=\sum_{j\in A} |g(v_{j})-g(v_{j-1})|+\sum_{j\in B}|g(v_{j})-g(v_{j-1})| \\
    &\leq s(1-\eta)+(r-s) = r-s\eta. 
\end{align*}
Consequently, 
\begin{equation}\label{eq:s.r-is-smaller-than-delta.eta}
    \frac{s}{r} < \frac{\delta}{\eta}.
\end{equation}
On the other hand, 
\begin{align*}
    \frac{|f(q)-f(p)|}{r} &\leq \sum_{j \in A} \frac{|f(v_j)-f(v_{j-1})|}{r} + \sum_{j \in B} \frac{|f(v_j)-f(v_{j-1 })|}{r} \\
    &\leq \frac{s}{r} \|f\|_L + \frac{r-s}{r} a^\eta_g (f) \leq \frac{\delta}{\eta} \|f\|_L + a^\eta_g (f),
\end{align*}
where the last inequality follows from \eqref{eq:s.r-is-smaller-than-delta.eta}. Since a pair of points $p\neq q \in V$ with $|g(q)-g(p)|>(1-\delta)d(p,q)$ was chosen arbitrarily, we conclude that $\nu_{L,g}^\delta (f) \leq a^\eta_g(f) + \tfrac{\delta}{\eta}\|f\|_L$.

(c) Define $\phi:\mathbb{R}\to\mathbb{R}$ by $\phi(t)=\text{dist}(t,\mathbb{Z})$. Then $\phi$ is $1$-Lipschitz and satisfies that 
\[
s,t\in\mathbb{R},\,|s-t|\leq 1 \implies |\phi(s)-\phi(t)|\leq 1-|s-t|.
\]
Let $h \in \Lip_0 (V,\mathbb{R})$ be defined as $h=\phi\circ g$. Then $\|h\|_L \leq 1$. If $uv \in E$, then $|g(v)-g(u)|\leq 1$; hence 
\[
|h(v)-h(u)|=|\phi(g(v))-\phi(g(u))|\leq 1-|g(v)-g(u)|.
\]
Part (a) therefore yields that $\|g\pm h \|_L \leq 1$. Since $g$ is an extreme point, we have $h=0$, i.e., $g(v) \in \mathbb{Z}$ for every $v \in V$.
\end{proof}

\begin{remark}\label{rem:strongly-exposed-edge}
    Given a graph $\Gamma=(V,E)$, the elementary molecule $m_{u,v}$ is a strongly exposed point of $B_{\mathcal{F}(V)}$ for every $uv \in E$. In fact, let $uv \in E$ be given, and consider 
    \[
    \phi = \frac{1}{2} (\chi_{\{u\}}-\chi_{\{v\}}) - \frac{1}{2} (\chi_{\{u\}}-\chi_{\{v\}})(0).
    \]
    It is straightforward to check that $\phi \in S_{\Lip_0(V,\mathbb{R})}$, $\phi(m_{u,v})=1$, and 
    \[
    |\phi(m_{p,q})|\leq \frac{1}{2} \quad\text{for every}~ (p,q)\notin \{(u,v),(v,u)\}.
    \]
    Since $B_{\mathcal{F}(V)}=\overline{\text{conv}}(\{m_{p,q}:pq \in E\})$ (by Lemma \ref{lem:basic-lemma-on-graphs}), we conclude that $\phi$ strongly exposes $B_{\mathcal{F}(V)}$ at $m_{u,v}$.
\end{remark}

The following theorem proves part~(a) of Theorem~\ref{thm:intro-graphs}. Let us mention that, for infinite graphs, passing from almost norming pairs to exactly saturated edges is not automatic. It becomes possible because positivity of the index forces $g$ to be integer-valued. This leads to an exact formula in terms of the saturated subgraph and restricts the possible positive values of the index.

\begin{theorem}\label{thm:Lipschitz-numerical-index-on-graphs}
    Let $\Gamma=(V,E)$ be a graph and let $g\in S_{\Lip_0 (V,\mathbb{R})}$. Then 
    \[
    n_{L,g}(V,\mathbb{R})=\frac{1}{D_g}.
    \]
    Moreover, 
    \[
    n_{L,g}(V,\mathbb{R})\in \{0,1\} \cup \left\{ \frac{1}{2k} : k \in \mathbb{N} \right\}.
    \]
\end{theorem}

\begin{proof}
For the lower bound, the case $D_g =\infty$ is trivial. Assume that $D_g<\infty$ and fix $f \in \Lip_0 (V,\mathbb{R})$. If $uv \in E_g$, then for every $0<\delta<1$, we have 
    \[
    |f(u)-f(v)|\leq \nu_{L,g}^\delta (f),
    \]
    which implies that $|f(u)-f(v)|\leq \nu_{L,g}(f)$. If $uv \in E$ and $u=v_0,v_1,\ldots,v_r=v$ is the shortest path connecting $u$ and $v$ with $v_jv_{j-1}\in E_g$, then $r=d_{H_g}(u,v)\leq D_g$; so 
    \begin{equation}\label{eq:estimate-fu-fv}
    |f(u)-f(v)|\leq \sum_{j=1}^r |f(v_{j})-f(v_{j-1})|  \leq r\nu_{L,g} (f) \leq D_g \nu_{L,g} (f).
    \end{equation}
    Since \eqref{eq:estimate-fu-fv} holds for every $uv \in E$, by Lemma \ref{lem:basic-lemma-on-graphs}, we obtain that $\|f\|_L \leq D_g \nu_{L,g}(f).$ This proves that 
    \[
    n_{L,g}(V,\mathbb{R}) \geq \frac{1}{D_g}.
    \]

    Next, we claim that $n_{L,g}(V,\mathbb{R})\leq \tfrac{1}{D_g}$. If $n_{L,g}(V,\mathbb{R})=0$, then there is nothing to prove; so we assume that $n_{L,g}(V,\mathbb{R})>0$. By Corollary \ref{cor:norm-derivative}(c), $g$ is an extreme point of $B_{\Lip_0 (V,\mathbb{R})}$. Lemma \ref{lem:basic-lemma-on-graphs}(c) therefore gives $g(V)\subseteq \mathbb{Z}$. 

    \emph{Step 1}. By Theorem \ref{thm:n-L-g-metric}, we have that 
    \[
    n_{L,g}(V,\mathbb{R}) = \inf_{0<\delta<1} \frac{1}{M_g^\delta}. 
    \]
    In particular, for $0<\delta<1$, we have that 
    \begin{equation}\label{eq:rho-g-delta-estimate-with-respect-to-d-distance}
    \rho_{g,\delta } (p,q) \leq \left( \frac{1}{n_{L,g} (V,\mathbb{R})} \right) d(p,q)
    \end{equation}
    for every $p,q\in V$. This also shows that $\rho_{g,\delta}(p,q)<\infty$ for every $p,q \in V$.

    \emph{Step 2}. Take $\delta>0$ sufficiently small so that $\delta < n_{L,g}(V,\mathbb{R})$. Let $uv \in E$ be given. By \eqref{eq:rho-g-delta-estimate-with-respect-to-d-distance}, we have 
    \[
    \rho_{g,\delta}(u,v)\leq \frac{1}{n_{L,g}(V,\mathbb{R})}.
    \]
    Since every admissible chain has a positive integer length and $\rho_{g,\delta}(u,v)<\infty$, the infimum defining $\rho_{g,\delta}(u,v)$ is attained. Thus, we can take $u=x_0,x_1,\ldots, x_m=v$ with $(x_{j-1},x_j)\in\mathfrak{N}_\delta (g)$ such that 
    \[
    \sum_{j=1}^m d(x_{j-1} ,x_j) = \rho_{g,\delta} (u,v).
    \] 
    Note that 
    \[
    0\leq d(x_{j-1},x_j)-|g(x_{j-1})-g(x_j)|< \delta d(x_{j-1},x_j) \leq \delta \left( \frac{1}{n_{L,g}(V,\mathbb{R})} \right) <1.
    \]
    Since $d(x_{j-1},x_j)\in\mathbb{N}$, and $g(V) \subseteq \mathbb{Z}$, we conclude that \[
    |g(x_{j-1} ) -g(x_j)|=d(x_{j-1},x_j).
    \]
    This shows that every edge of a shortest path connecting $x_{j-1}$ and $x_j$ is $g$-saturated; hence 
    \[
    d_{H_g}(u,v) \leq \sum_{j=1}^m d(x_{j-1}, x_j) = \rho_{g,\delta} (u,v) \leq \frac{1}{n_{L,g}(V,\mathbb{R})}.
    \]
    Since $uv\in E$ was chosen arbitrarily, we conclude that 
    \[
    D_g \leq \frac{1}{n_{L,g}(V,\mathbb{R})}.
    \]

    To prove the ``moreover'' part, suppose that $n_{L,g}(V,\mathbb{R})>0$. Since $g(V)\subseteq \mathbb{Z}$ in this case, we have $g(v)-g(u) \in \{-1,0,1\}$ for every $uv \in E$. Also, by Step 2 above, we have
    $$d_{H_g}(u,v) \leq \frac{1}{n_{L,g}(V,\mathbb{R})}.$$
    If $d_{H_g} (u,v)=2k+1$ for some $k\in\mathbb{N}$, then there exist $u=x_0,x_1,\ldots,x_{2k+1}=v$ such that $x_jx_{j-1}\in E_g$ for every $j$. Since $g(x_j)-g(x_{j-1})\in\{-1,1\}$ for $j=1,\ldots, 2k+1$, 
    \begin{align*}
    g(v)-g(u)&=\sum_{j=1}^{2k+1} (g(x_j)-g(x_{j-1})) \\
    &\in \{-2k-1,-2k+1,\ldots,-1,1,\ldots,2k-1,2k+1\}.
    \end{align*}
    On the other hand, since $|g(v)-g(u)|\leq d(u,v)=1$, we conclude that $g(v)-g(u)$ must be $1$ or $-1$. This implies that $uv \in E_g$ and $d_{H_g}(u,v)=1$, which is a contradiction. It follows that $d_{H_g}(u,v) \in \{1\}\cup 2\mathbb{N}$, and that $D_g \in \{1\}\cup 2 \mathbb{N}$.
\end{proof}

\begin{remark}
The restriction on the possible values in Theorem~\ref{thm:Lipschitz-numerical-index-on-graphs} does not extend to arbitrary finite metric spaces. Indeed, for each $\tau\in(0,1]$, Example~\ref{ex:three-point-example} provides a three-point metric space $M_\tau$ and $g\in S_{\Lip_0(M_\tau,\mathbb R)}$ such that $n_{L,g}(M_\tau,\mathbb R)=\tau$.
\end{remark}

\subsection{Extreme and strongly extreme points}

The saturated subgraph also describes the extremal structure of the unit ball of $\Lip_0(V,\mathbb{R})$. Extremality requires $H_g$ to be connected, whereas strong extremality requires a uniform bound on the distances in $H_g$ between endpoints of edges of $\Gamma$. 

    \begin{theorem}\label{thm:characterization-of-extreme-and-strongly-extreme-function-in-Lip-spaces}
    Let $\Gamma=(V,E)$ be a graph and let $g\in S_{\Lip_0 (V,\mathbb{R})}$. Then 
        \begin{enumerate}
            \item[\textup{(a)}] $g$ is an extreme point of  $B_{\Lip_0(V,\mathbb{R})}$ if and only if the subgraph $H_g$ is connected. In this case, $g(V)\subseteq \mathbb{Z}$.
            \item[\textup{(b)}] $g$ is a strongly extreme point of $B_{\Lip_0(V,\mathbb{R})}$ if and only if $n_{L,g}(V,\mathbb{R})>0$. In this case,
            \begin{equation}\label{eq:g-pm-h-estimate-wrt-Dg}
            \max \{\|g+h\|_L,\|g-h\|_L\}\geq 1+ \frac{\|h\|_L}{D_g} \quad\text{for every}~ h\in\Lip_0 (V,\mathbb{R}).
            \end{equation}
        \end{enumerate}
    \end{theorem}

\begin{proof}
(a) Since $\max \{|a+b|,|a-b|\}=|a|+|b|$ for all $a,b \in \mathbb{R}$, we obtain from Lemma \ref{lem:basic-lemma-on-graphs} that 
\begin{equation}\label{eq:g-pm-h-estimate}
    \max \{ \|g+h\|_L, \|g-h\|_L\} = \sup_{uv\in E} (|g(v)-g(u)|+|h(v)-h(u)|).
\end{equation}
We first claim that if $H_g$ is connected, then $g$ is an extreme point. Suppose that $g\pm h\in B_{\Lip_0 (V,\mathbb{R})}$ for some $h \in \Lip_0(V,\mathbb{R})$. If $uv \in E_g$, i.e., $|g(v)-g(u)|=1$, then by \eqref{eq:g-pm-h-estimate}
\[
1+|h(v)-h(u)| \leq \max \{\|g+h\|_L,\|g-h\|_L\} \leq 1;
\]
hence $h(u)=h(v)$. Since $H_g$ is connected, $h(0)=0$, and $h$ is constant on $E_g$, we conclude that $h =0$. It follows that $g$ is an extreme point. 

Conversely, assume that $g$ is an extreme point. By Lemma \ref{lem:basic-lemma-on-graphs}(c), $g(V)\subseteq \mathbb{Z}$.

If $H_g$ is not connected, then there exists a connected component $C$ that does not contain $0$. Consider $h= (\tfrac12)\chi_C \in \Lip_0 (V,\mathbb{R})$. Given $uv \in E$, if $u,v\in C$ or $u,v\not\in C$, then $h(v)-h(u)=0$. If $u \in C$ and $v \not\in C$, i.e., $u$ and $v$ are not connected by edges in $E_g$, then $g(u)=g(v)$ since $g(u)-g(v)\in\{-1,0,1\}.$ It follows from \eqref{eq:g-pm-h-estimate} that 
\[
\|g\pm h\|_L = \sup_{uv\in E} (|g(v)-g(u)|+|h(v)-h(u)|) \leq 1.
\]
This contradicts the extremality of $g$; so $H_g$ is connected.

(b) By Corollary \ref{cor:norm-derivative}, it is enough to prove that if $g$ is a strongly extreme point, then $n_{L,g}(V,\mathbb{R})>0$. Since $g$ is, in particular, an extreme point, by the above argument, we have that $H_g$ is connected and $g(v)-g(u)\in\{-1,0,1\}$ for every $uv \in E$. Assume to the contrary that $D_g = \infty$. For each $n \in \mathbb{N}$, find $u_nv_n \in E$ such that $L_n:=d_{H_g}(u_n,v_n)>n$. Consider 
\[
h_n(z):= \frac{\min \{ d_{H_g}(z,v_n),L_n \}-\min \{d_{H_g}(0,v_n),L_n\}}{L_n}
\]
for each $z\in V$. Note that $\|h_n\|_L=1$ and $|h_n(u_n)-h_n(v_n)|=1=d(u_n,v_n)$. If $xy \in E_g$, then 
\[
|d_{H_g} (x,v_n)-d_{H_g}(y,v_n)|\leq d_{H_g} (x,y)\leq 1;
\]
hence 
\[
|h_n(y)-h_n(x)|\leq \frac{1}{L_n}.
\]
If $xy\in E\setminus E_g$, then $g(x)=g(y)$, since $g(x)-g(y) \in \{-1,0,1\}$. It follows that 
\[
\|g + h_n \|_L, \|g-h_n \|_L \leq 1+ \frac{1}{L_n}.
\]
Define 
\[
f_n:=\frac{g+h_n}{1+1/L_n} \quad \text{and} \quad k_n:=\frac{g-h_n}{1+1/L_n}.
\]
Then $f_n, k_n \in B_{\Lip_0(V,\mathbb{R})}$, and 
\[
\frac{f_n+k_n}{2} = \frac{g}{1+1/L_n} \to g.
\]
However, 
\[
\|f_n-k_n\|_L = \frac{2\|h_n\|_L}{1+1/L_n}=\frac{2}{1+1/L_n}\to 2.
\]
This contradicts that $g$ is a strongly extreme point.

Finally, we prove the estimate \eqref{eq:g-pm-h-estimate-wrt-Dg}. Let $uv \in E$ be given. Since $D_g<\infty$, there exists $u=v_0,v_1,\ldots,v_r=v$ with $v_jv_{j-1}\in E_g$ such that 
\[
d_{H_g}(u,v)=r \leq D_g.
\]
Notice that for each $j=1,\ldots,r$, we have from \eqref{eq:g-pm-h-estimate} that 
\[
1+|h(v_j)-h(v_{j-1})| \leq \max \{\|g+h\|_L,\|g-h\|_L\}.
\]
It follows that 
\begin{align*}
|h(v)-h(u)|&\leq \sum_{j=1}^r |h(v_j)-h(v_{j-1})| \\
&\leq d_{H_g}(u,v) (\max \{\|g+h\|_L,\|g-h\|_L\}-1) \\
&\leq D_g (\max \{\|g+h\|_L,\|g-h\|_L\}-1) .
\end{align*}
Taking the supremum over $uv\in E$, we obtain 
\[
\|h\|_L\leq D_g(\max \{\|g+h\|_L,\|g-h\|_L\}-1),
\]
which is equivalent to \eqref{eq:g-pm-h-estimate-wrt-Dg}.
\end{proof}

\subsection{Bipartite graphs}

A \emph{bipartite graph} is a graph whose vertex set can be partitioned as $V=A \sqcup B$ so that every edge joins a vertex of $A$ to a vertex of $B$.

\begin{theorem}\label{thm:bipartite-characterization}
    Given a connected graph $\Gamma=(V,E)$, $\Gamma$ is bipartite if and only if there exists $g \in S_{\Lip_0(V,\mathbb{R})}$ such that $n_{L,g}(V,\mathbb{R})=1$. Moreover, if $\Gamma$ is a bipartite graph, then, for every $g \in S_{\Lip_0(V,\mathbb{R})}$, the following are equivalent:
    \begin{enumerate}
        \item[\textup{(a)}] $g$ is an extreme point of $B_{\Lip_0(V,\mathbb{R})}$; 
        \item[\textup{(b)}] $g$ is a strongly extreme point of $B_{\Lip_0(V,\mathbb{R})}$; 
        \item[\textup{(c)}] $n_{L,g}(V,\mathbb{R})>0$;
        \item[\textup{(d)}] $n_{L,g}(V,\mathbb{R})=1$;
        \item[\textup{(e)}] $|g(u)-g(v)|=1$ for every $uv \in E$.
    \end{enumerate}
    In particular, $n_{L,g}(V,\mathbb{R})\in\{0,1\}$ and every extreme point of $B_{\Lip_0(V,\mathbb{R})}$ is a spear element.
\end{theorem}

\begin{proof}
Suppose first that $\Gamma$ is bipartite, with $V=A\sqcup B$ and $0\in A$. For $g=\chi_B$, every edge is $g$-saturated, so $D_g=1$ and Theorem~\ref{thm:Lipschitz-numerical-index-on-graphs} gives $n_{L,g}(V,\mathbb R)=1$. 

Conversely, if $n_{L,g}(V,\mathbb R)=1$, then Theorem~\ref{thm:Lipschitz-numerical-index-on-graphs} gives $D_g=1$, and hence every edge is $g$-saturated. Since $g(0)=0$ and $\Gamma$ is connected, $g(V)\subseteq\mathbb Z$. Thus $g^{-1}(2\mathbb Z)$ and $g^{-1}(2\mathbb Z+1)$ form a bipartition of $\Gamma$.

Next, suppose that $\Gamma$ is a bipartite graph. Then 
    (d) $\implies$ (c), and (b) $\implies$ (a) are clear. (c) $\implies$ (b) is a consequence of Theorem \ref{thm:characterization-of-extreme-and-strongly-extreme-function-in-Lip-spaces}. If (e) holds, then $E=E_g$; hence $D_g=1$. Thus, Theorem \ref{thm:Lipschitz-numerical-index-on-graphs} yields that $n_{L,g}(V,\mathbb{R})=1$. This proves (e) $\implies$ (d). It remains to prove (a) $\implies$ (e). Suppose that $g$ is an extreme point. Then Theorem \ref{thm:characterization-of-extreme-and-strongly-extreme-function-in-Lip-spaces} implies that $H_g$ is connected. Let $uv \in E$ be given. It is enough to check that $uv \in E_g$. Let $u=v_0,v_1,\ldots, v_r=v$ be the shortest path in $H_g$. Since every edge in $\Gamma$ connects points that are in different parts of the bipartition, $r$ is an odd number. It follows that 
    \[
    g(v)-g(u)=\sum_{j=1}^r (g(v_j)-g(v_{j-1}))
    \]
    is a sum of an odd number of terms, each equal to $1$ or $-1$. Therefore $g(v)-g(u)$ is an odd integer. However, since $uv \in E$ (so, $|g(v)-g(u)|\leq 1$), we conclude that $g(v)-g(u)=\pm 1$. This proves that $uv \in E_g$. 
\end{proof}

\subsection{Classical numerical indices}

We now pass from the index associated with a fixed Lipschitz function to the classical numerical indices of $\Lip_0 (V,\mathbb{R})$ and $\mathcal{F}(V)$. The main tool is a rounding argument: the unit ball of $\Lip_0(V,\mathbb{R})$ is the norm-closed convex hull of its integer-valued elements.

The next two theorems prove part~(b) of Theorem~\ref{thm:intro-graphs}. The first gives the bounds for the classical numerical indices, and the second characterizes bipartiteness in terms of these indices. 

\begin{theorem}\label{thm:general-lower-bound-for-numerical-index-of-LF-graph}
    Let $\Gamma=(V,E)$ be a connected graph. Then 
    \[
    \frac{1}{2} \leq n(\Lip_0 (V,\mathbb{R}))\leq n(\mathcal{F}(V)) \leq 1.
    \]
\end{theorem}

\begin{proof}
    First, we claim that every element $f$ in $\Lip_0(V,\mathbb{R})$ can be approximated by convex combinations of integer-valued elements in $\Lip_0(V,\mathbb{R})$ in norm. As a matter of fact, assume that $\|f\|_L=1$, and for $0\leq t <1$, consider 
    \[
    q_t (v) := \lfloor f(v)+t\rfloor
    \]
    for each $v\in V$. Then for every $uv\in E$, since $|f(u)-f(v)|\le 1$, we have 
    \[
    |q_t (u)-q_t (v)|\le 1;
    \]
    hence $q_t$ is also $1$-Lipschitz, and $q_t(V)\subseteq \mathbb{Z}$. Moreover, consider 
    \[
    Q_N := \frac{1}{N}\sum_{j=0}^{N-1} q_{\tfrac{j}{N}} \in \conv (\{ g\in B_{\Lip_0 (V,\mathbb{R})} : g(V)\subseteq \mathbb{Z} \} ).
    \]
Observe that 
\[
Q_N (v) = \frac{\lfloor Nf(v)\rfloor}{N};
\]
hence 
\[
\|Q_N-f\|_L \leq \frac{1}{N}.
\]
It follows that $f \in \overline{\conv} (\{ g\in B_{\Lip_0 (V,\mathbb{R})} : g(V)\subseteq \mathbb{Z} \} )$. Since $f$ was arbitrarily chosen, we conclude that
\begin{equation}\label{eq:integer-valued-convex-hull}
    B_{\Lip_0 (V,\mathbb{R})}=\overline{\conv} (\{ g\in B_{\Lip_0 (V,\mathbb{R})} : g(V)\subseteq \mathbb{Z} \} ).
\end{equation}

Second, let $T\in \mathcal{L}(\Lip_0(V,\R))$ be given. Let $g \in B_{\Lip_0 (V,\mathbb{R})}$ be given such that $g(V)\subseteq \mathbb{Z}$. Let $uv \in E$. Then $g(u)-g(v)\in \{-1,0,1\}$. Consider $\phi \in S_{\Lip_0(V,\mathbb{R})}$ from Remark \ref{rem:strongly-exposed-edge} that satisfies 
\[
\phi(m_{u,v})=1 \quad \text{and} \quad |\phi(m_{p,q})|\leq \frac{1}{2}, \quad \forall (p,q)\not\in \{(u,v),(v,u)\}.
\]
If $|g(u)-g(v)|=1$, then we have either $g(\pm m_{u,v}) = 1$. Thus, we have  
\begin{equation}\label{eq:nu(T)-first-estimate}
\nu(T) \ge |(Tg) (m_{u,v})|= |(Tg)(u)-(Tg)(v)|.
\end{equation}
If $g(u)-g(v)=0$, then notice that $\frac{1}{2}g \pm  \phi$ satisfies that 
\[
\left(\frac{1}{2}g \pm  \phi\right)(m_{u,v}) = \pm 1,  \quad \text{and} \quad \left|\left(\frac{1}{2}g \pm  \phi\right)(m_{p,q})\right| \leq \frac{1}{2}+\frac{1}{2}=1
\]
for all $(p,q)\notin \{(u,v),(v,u)\}$.
It follows that 
\begin{equation}\label{eq:nu(T)-second-estimate}
\nu(T) \geq \max_{\pm} \left| \left\langle T\left(\frac{1}{2}g \pm  \phi\right) , m_{u,v} \right\rangle \right| \ge \frac{1}{2} |(Tg)(m_{u,v})|.
\end{equation}
Combining \eqref{eq:nu(T)-first-estimate} and \eqref{eq:nu(T)-second-estimate}, we conclude that 
\[
\nu(T) \geq \frac{1}{2} |(Tg)(u)-(Tg)(v)|, \quad \forall uv\in E.
\]
Hence $\nu(T)\geq \tfrac{\|Tg\|}{2}$ for every integer-valued
$g\in B_{\operatorname{Lip}_0(V,\mathbb R)}$. The first claim \eqref{eq:integer-valued-convex-hull} guarantees that $n(\Lip_0(V, \mathbb{R}))\geq \tfrac12$.

The second inequality in the statement follows from a general fact for Banach spaces.
\end{proof}

\begin{theorem}\label{thm:numerical-index-of-LF-bipartite-graphs} 
   Let $\Gamma=(V,E)$ be a connected graph. Then $\Gamma$ is a bipartite graph if and only if $n(\mathcal{F}(V))=1$ if and only if $n(\Lip_0(V,\mathbb{R}))=1$. 
\end{theorem}

\begin{proof}
Suppose that $\Gamma$ is a bipartite graph. Let $T : \Lip_0 (V,\mathbb{R})\to\Lip_0(V,\mathbb{R})$ be a bounded linear operator. Let $g$ be an extreme point of $B_{\Lip_0 (V,\mathbb{R})}$. For any $uv \in E$, by Theorem \ref{thm:bipartite-characterization}, $\phi = \delta_v-\delta_u$ or $\delta_u-\delta_v$ satisfies that $\phi(g)=1$. Then we have 
\[
\nu(T)\geq |\phi(Tg)| = |(Tg)(u)-(Tg)(v)|.
\]
Since it holds for every $uv \in E$, we conclude that $\nu(T)\geq \|Tg\|$. 

Next, we claim that $B_{\Lip_0 (V,\mathbb{R})}=\overline{\text{conv}}(\text{ext}(B_{\Lip_0 (V,\mathbb{R})}))$, which would imply that $\nu(T) \geq \|Tf\|$ for every $f \in B_{\Lip_0 (V,\mathbb{R})}$, equivalently, $n( \Lip_0(V,\mathbb{R}))=1$. Let $f \in B_{\Lip_0 (V,\mathbb{R})}$ be given. Write $V = A \sqcup B$ with $0 \in A$ and let 
\[
\sigma(v)=\begin{cases}
    0, \quad &v\in A, \\
    1, \quad &v \in B.
\end{cases}
\]
For $0\leq t < 1$, define $h_t \in \Lip_0(V,\mathbb{R})$ by  
\[
h_t (v) = 2 \left \lfloor \frac{f(v)+\sigma(v)}{2}+t \right\rfloor-\sigma (v).
\]
If $u\in A$, $v \in B$ and $uv \in E$, then since $|f(u)-f(v)|\leq 1$, we have 
\[
h_t (v) - h_t(u) = \left( 2 \left \lfloor \frac{f(v)+1}{2}+t \right\rfloor -1 \right) - 2\left \lfloor \frac{f(u)}{2}+t \right\rfloor \in \{-1,1\}.
\]
It follows from Theorem \ref{thm:bipartite-characterization} that $h_t$ is an extreme point of $B_{\Lip_0 (V,\mathbb{R})}$. 

Moreover, if $$
Q_N = \frac{1}{N}\sum_{j=0}^{N-1}h_{\tfrac{j}{N}} \in \text{conv}(\text{ext}(B_{\Lip_0 (V,\mathbb{R})})),
$$ then we can rewrite $Q_N$ as 
\[
Q_N (v) = \frac{2}{N} \left\lfloor \frac{N(f(v)+\sigma(v))}{2} \right\rfloor -\sigma(v)
\]
for every $v\in V$. It follows that 
\[
-\frac{2}{N}<Q_N(v)-f(v)\leq 0;
\]
hence Lemma \ref{lem:basic-lemma-on-graphs} yields that $\|Q_N-f\|_L \to 0$.

Finally, suppose that $n(\mathcal{F}(V))=1$. Then for every extreme point $\phi \in B_{\Lip_0(V,\mathbb{R})}$ and for every denting point $\mu$ of $B_{\mathcal{F}(V)}$, we have $|\phi(\mu)|=1$ \cite[Lemma 1.(ii)]{LMP1999}. Let $g$ be an extreme point of $B_{\Lip_0(V,\mathbb{R})}$. By Remark \ref{rem:strongly-exposed-edge}, $m_{u,v}$ is a strongly exposed point for every $uv \in E$. Hence, we obtain that $|g(m_{u,v})|=1$ for all $uv \in E$. Since $g(0)=0$ and $\Gamma$ is connected, we have then $g(V) \subseteq \mathbb{Z}$. Thus, $V$ is decomposed into $g^{-1}(2\mathbb{Z})$ and $g^{-1}(2\mathbb{Z}+1)$, and this proves that $\Gamma$ is bipartite.
\end{proof}

\subsection{Cycle graphs}

A \emph{cycle} is a closed path that returns to its starting point at the end, and otherwise has no repeated vertices. A \emph{cycle graph} is a graph that consists of a single cycle. The cycle graph with $m$ edges is denoted by $C_m$. 

\begin{theorem}\label{thm:cycle-graph-characterization}
    Let $C_m=(\{v_0,\ldots,v_m\},E)$ be a cycle graph with $v_m=v_0$, and let $g \in S_{\Lip_0(C_m,\mathbb{R})}$. For $j=0,\ldots, m-1$, put $a_j=g(v_{j+1})-g(v_{j})$. 
    \begin{enumerate}
        \item[\textup{(a)}] If $m=2k$ for some $k\geq 2$, then 
        \[
        n_{L,g}(C_{2k},\mathbb{R})=\begin{cases}
            1, \quad &\text{$|a_j|=1$ for every $j$,} \\ 
            0, \quad &\text{otherwise.}
        \end{cases}
        \]
        \item[\textup{(b)}] If $m=2k+1$ for some $k\geq 1$, then 
        \[
        n_{L,g}(C_{2k+1},\mathbb{R})=\begin{cases}
            \dfrac{1}{2k}, \quad &\text{exactly one $a_j=0$ and all other $a_i \in \{-1,1\}$,} \\ 
            0, \quad &\text{otherwise.}
        \end{cases}
        \]
    \end{enumerate}
\end{theorem}

\begin{proof}
    Let $f \in \Lip_0(C_m,\mathbb{R})$ be given. Let $b_j:=f(v_{j+1})-f(v_j)$. Note that $\|f\|_L=\max_{0\le j \le m-1} |b_j|$. If we let 
    \[
    S=\{j: |a_j|=1\} \quad \text{and} \quad U=\{j:|a_j|<1\},
    \]
    then $\nu_{L,g}(f)=\max_{j \in S} |b_j|$. It follows that 
    \[
    n_{L,g}(C_m,\mathbb{R})=\inf\left \{\max_{j\in S} |b_j| : \max_{0\le j \le m-1} |b_j|=1, \, \sum_{j=0}^{m-1} b_j=0 \right\} .
    \]
    We claim that 
    \[
    n_{L,g}(C_{m},\mathbb{R})= \begin{cases}
            1, \quad &\text{$|U|=0$}, \\ 
            \dfrac{1}{m-1}, \quad &\text{$|U|=1$}, \\
            0, \quad &\text{$|U|\ge 2$}. \\ 
        \end{cases}
    \]
    
    If $|U|=0$, i.e., every edge is $g$-saturated; then $\|f\|_L=\nu_{L,g}(f)$ for every $f \in \Lip_0(C_m,\mathbb{R})$. This shows that $n_{L,g}(C_m,\mathbb{R})=1$. 

    Suppose that $|U|=1$. For simplicity, let $v_{m-1}v_0$ be the only unsaturated edge. Then $d_{H_g}(v_{m-1},v_0)=m-1$, as $v_0$ and $v_{m-1}$ can only be joined by $m-1$ edges $v_0v_1, \ldots, v_{m-2}v_{m-1}$. For every edge different from $v_{m-1}v_0$, $d_{H_g}$ has the value $1$. It follows that $D_g=\max \{1,m-1\}=m-1$; hence Theorem \ref{thm:Lipschitz-numerical-index-on-graphs} shows that $n_{L,g}(C_m,\mathbb{R})=\tfrac{1}{(m-1)}$.

    Finally, suppose that $|U|\ge 2$, i.e., $C_m$ has at least two unsaturated edges. Then the graph $H_g$ is disconnected since for any unsaturated edge $uv$, the only alternative path in the cycle joining $u$ and $v$ contains another unsaturated edge, hence $u$ and $v$ lie in different connected components of $H_g$. It follows that $d_{H_g}(u,v) = \infty$; hence $D_g=\infty$. So, Theorem \ref{thm:Lipschitz-numerical-index-on-graphs} shows that $n_{L,g}(C_m,\mathbb{R})=0$.

    We finish the proof by mentioning that if $|U|=0$, then since $\sum_{j=0}^{m-1} a_j=0$ with $a_j\in\{-1,1\}$, the index $m$ must be even. If $|U|=1$, then $m$ is an odd number, and there exists $r\in \{0,\ldots,m-1\}$ such that 
    \[
    |a_r|<1 \quad \text{and} \quad a_r = -\sum_{j\neq r} a_j \,\,\text{ with $a_j \in \{-1,1\}$, $\forall j\neq r$}. 
    \]
    This shows that $a_r \in \mathbb{Z}$, which forces $a_r=0$. 
    \end{proof}

\begin{theorem}\label{thm:numerical-index-of-LF-cycle-graphs}
Let $C_m=(\{v_0,\ldots,v_m\},E)$ be a cycle graph with $v_m=v_0$. For $m \geq 3$, 
\[\
n(\mathcal{F}(C_m))= n(\Lip_0 (C_m,\mathbb{R})) = \begin{cases}
    1, \quad &\text{$m$ is even},\\
    \frac{m-2}{m-1}, \quad &\text{$m$ is odd.}
\end{cases}
\]
\end{theorem}

\begin{proof}
If $m$ is an even number, then $C_m$ is a bipartite graph. Thus, Theorem \ref{thm:numerical-index-of-LF-bipartite-graphs} yields that $n(\mathcal{F}(C_m))=n(\Lip_0(C_m,\mathbb{R}))=1$.

    Suppose that $m=2k+1$ for some $k \in \mathbb{N}$. Note that $\sum_{j=1}^{m} (f(v_{j})-f(v_{j-1}) ) = -f(v_0)+f(v_m)=0$. Moreover, the map
    \begin{align*}
    \Lip_0(C_m,\mathbb{R}) &\to E_m:=\left \{x=(x_1,\ldots,x_m)\in\mathbb{R}^m:\sum_{j=1}^m x_j=0\right \}, \\
    &f\mapsto (f(v_j)-f(v_{j-1}))_{j=1}^m
    \end{align*}
    is a surjective linear isometry, where $E_m$ is endowed with $\|\cdot\|_\infty$-norm. Thus, it is enough to prove that $n(E_m)=\frac{m-2}{m-1}$. 

    \emph{Step 1.} It is routine to observe that 
    \begin{align*}
    &\text{ext}(B_{E_m}) \\
    &=\{(x_i)_{i=1}^m\in E_m : |x_i|=1 \text{ for $i\in I$ with $|I|=2k$ and $x_j=0$ for $j\notin I$}\}.
    \end{align*}
    Obviously, for $x=(x_i)_{i=1}^m \in \text{ext}(B_{E_m})$, we have 
    \[
    |\{i\in \{1,\ldots,m\}:x_i=1\}|=|\{i\in \{1,\ldots,m\}:x_i=-1\}|=k.
    \]
Since $B_{E_m^*} = \conv \{\pm e_i^*\restricted_{E_m} : 1 \le i \le m\}$, the numerical radius can be computed from the norming coordinate functionals. That is, for a bounded linear operator $T:E_m \to E_m$, we have 
    \begin{equation}\label{eq:nu-T-description}
    \nu(T) = \max \{| (Tx)_i|: x \in B_{E_m}, \,1\le i \le m, \,|x_i|=1\}.
    \end{equation}

    \emph{Step 2.} We claim that $n(E_m)\geq \frac{m-2}{m-1}$. Let $T \in \mathcal{L}(E_m,E_m)$ be given and $x \in \text{ext}(B_{E_m})$ such that $\|Tx\|=\|T\|$. Let $j \in \{1,\ldots,m\}$ such that $|(Tx)_j|=\|T\|$.
    \begin{itemize}
        \item If $|x_j|=1$, then we have $\nu(T)=\|T\|$ by \eqref{eq:nu-T-description}.
        \item If $x_j=0$, then define $z=(z_n)_{n=1}^m \in E_m$ by 
        \[
        z_j =1 \quad \text{and} \quad z_i = -\frac{1}{m-1} \,\, \text{ for $i \neq j$}.
        \]
        Consider 
        \[
        u=\left(\frac{m-2}{m-1} \right)x + z \quad \text{and} \quad v=\left(\frac{m-2}{m-1} \right)x -z.
        \]
        Then $u,v \in B_{E_m}$ and $u_j=-v_j=1$. Note from \eqref{eq:nu-T-description} that 
        \begin{align*}
        &\nu(T) \geq |(Tu)_j|=\left| \left(\frac{m-2}{m-1} \right) (Tx)_j + (Tz)_j \right|,\quad\text{and}  \\
        &\nu(T)\geq |(Tv)_j|= \left| \left(\frac{m-2}{m-1} \right) (Tx)_j - (Tz)_j \right|.
        \end{align*}
    \end{itemize}
    It follows that 
    \[
    \nu(T)\geq \left(\frac{m-2}{m-1}\right)|(Tx)_j| = \left(\frac{m-2}{m-1}\right) \|T\|.
    \]

    \emph{Step 3.} Finally, we claim that $n(E_m)\leq \frac{m-2}{m-1}.$ Consider $T\in\mathcal{L}(E_m,E_m)$ defined by 
    \begin{align*}
    Tx=x^*(x)z, \,\, \text{where} \,\, x^*(x)&=\frac{1}{2k}\left(\sum_{i=1}^k x_i -\sum_{i=k+1}^{2k} x_i \right),\quad\text{and} \\
    &z = \left(-\frac{1}{2k},\ldots,-\frac{1}{2k}, 1\right) \in S_{E_m}.
    \end{align*}
    Then $T$ is a rank-one norm-one operator, and one can check that $|(Tx)_i|\leq \frac{m-2}{m-1}$ for every $x \in B_{E_m}$ with $|x_i|=1$. Thus, \eqref{eq:nu-T-description} yields that $\nu(T)\leq \frac{m-2}{m-1}.$
\end{proof}

Theorem \ref{thm:cycle-graph-characterization} and Theorem \ref{thm:numerical-index-of-LF-cycle-graphs} highlight the difference between the behaviours of $n_{L,g}(\cdot)$ and $n(\mathcal{F}(\cdot))$: for cycle graphs $C_{2k+1}$, we have 
\[
\sup_{g \in S_{\Lip_0(C_{2k+1},\mathbb{R})}} n_{L,g}(C_{2k+1},\mathbb{R}) = \frac{1}{2k}\to0~\text{while}~n(\mathcal{F}(C_{2k+1}))=\frac{2k-1}{2k}\to1.
\]

\section{Attainment of the Lipschitz numerical radius}\label{sec:radius-attaining}

Throughout this section, $Y$ denotes a \emph{real or complex}
Banach space. For $g\in S_{\Lipz(M,Y)}$, consider the following set of $g$-norming triples: 
\begin{align*}
\Pi_g &= \bigcap_{0<\delta<1} \Pi_g (\delta) \\
&= \left \{ (p,q,y^*) \in M\times M \times S_{Y^*}: p\neq q, \, y^*\left(\frac{g(p)-g(q)}{d(p,q)} \right )=1 \right \}. 
\end{align*}
Note that $\Pi_g \neq \emptyset $ if and only if there exist $p\neq q \in M$ such that $\|g(p)-g(q)\|=d(p,q)$, i.e., $g$ \emph{strongly attains its Lipschitz norm}.

\begin{definition}\label{def:g-Lipschitz-numerical-radius-attaining}
    Let $M$ be a metric space, $Y$ a Banach space, and $g \in \Lip_0(M,Y)$ satisfy $\|g\|_L=1$. A map $f \in \Lip_0 (M,Y)$ is said to \emph{attain its $g$-Lipschitz numerical radius} if there exists $(p,q,y^*)\in \Pi_g$ such that 
    \[
    \left|y^* \left( \frac{f(p)-f(q)}{d(p,q)} \right) \right|=\nu_{L,g}(f).
    \]
    We denote by $\LipNRA_g (M,Y)$ the set of all Lipschitz maps from $M$ to $Y$ that attain their $g$-Lipschitz numerical radii. 
\end{definition}

\subsection{Finiteness and universal attainment}

We first show that every infinite metric space $M$ admits a strongly norm-attaining reference map $g \in S_{\Lip_0(M,\mathbb{K})}$ for which some Lipschitz function fails to attain its $g$-Lipschitz numerical radius.

\begin{proposition}\label{prop:non-isolated-point}
    If $M$ is an infinite metric space, then 
    there exist $f,g \in S_{\Lip_0 (M,\mathbb{K})}$ such that 
    \begin{itemize}
        \item $g$ strongly attains its Lipschitz norm, 
        \item $f$ does not strongly attain its norm, and 
        \item $\nu_{L,g}(f)=1$.
    \end{itemize}
    In particular, $f$ does not attain its $g$-Lipschitz numerical radius. 
\end{proposition}

\begin{proof} We construct real-valued functions, which we regard as $\mathbb K$-valued when necessary. We divide the proof into two cases.

Case 1: $M$ contains a non-isolated point. Let $p_0 \in M$ be an accumulation point of $M$. Let $(p_n)$ be a sequence in $M\setminus \{p_0\}$ converging to $p_0$. Define $g \in \Lip_0(M,\mathbb{R})$ by $g(p)=d(p,p_0)-d(0,p_0)$ for every $p \in M$. Then $\|g\|_L=1$ and, by definition, $g(p)-g(p_0)=d(p,p_0)$ for every $p \in M$. So, $g$ strongly attains its Lipschitz norm. 

Define $f \in \Lip_0 (M,\mathbb{R})$ by 
    \[
    f(p) = \sin \left( d(p,p_0) \right) - \sin (d(0,p_0))
    \]
    for $p\in M$. Note that $\|f\|_L\leq 1$. Moreover, since $d(p_n,p_0)\to0$, we have 
    \begin{equation}\label{eq:sine-slope-tend-to-1}
    \frac{f(p_n)-f(p_0)}{d(p_n,p_0)} = \frac{\sin (d(p_n,p_0))}{d(p_n,p_0)} \to 1.
    \end{equation}
    So, $\|f\|_L=1$. Since each $(p_n,p_0)$ is a pair where $g$ attains its slope, we conclude from \eqref{eq:sine-slope-tend-to-1} that $\nu_{L,g}(f) \geq 1$. Since $\nu_{L,g}(f)\leq \|f\|_L \leq 1$, we have $\nu_{L,g}(f)=1.$ Finally, since $|\sin (x)-\sin(y)|<|x-y|$ for every $x\neq y$, we have that 
    \[
    |f(p)-f(q)|<|d(p,p_0)-d(q,p_0)| \leq d(p,q)
    \]
    for $p\neq q \in M$ with $d(p,p_0)\ne d(q,p_0)$. 
    This finishes the proof in this case.

    Case 2: Every point of $M$ is an isolated point. Since $M$ is infinite, $\mathcal F(M)$ contains a complemented copy of $\ell_1$ (see \cite{CDW}). In particular, $\mathcal F(M)$ is non-reflexive; hence there exists $\phi\in S_{\Lip_0(M,\mathbb R)}$ which does not strongly attain its norm. Let $p_0 \in M\setminus \{0\}$ be fixed, and consider 
    \[
    r:= \inf_{q \ne p_0 } d(p_0, q) >0 \, \text{and} \, s:=\sup_{q\neq p_0} \frac{|\phi(p_0)-\phi(q)|}{d(p_0,q)}.
    \]
    If $s=1$, then we can take $f:=\phi$ and $g:=d(\cdot ,p_0)-d(0,p_0)$ to verify the assertion. Suppose that $s<1$. Define 
    \[
    h(x) =
    \begin{cases}
        \phi(x) \quad &\text{if } x \ne p_0, \\ 
        \inf_{q\neq p_0} (\phi(q)+d(q,p_0) ) =:A \quad &\text{if } x = p_0.
    \end{cases}
    \]
    Observe that for every $p\in M\setminus \{p_0\}$, 
    \begin{align*} 
        \phi(p)+d(p,p_0)-\phi(p_0) 
        &=d(p,p_0)-(\phi(p_0)-\phi(p)) \\
        &\ge d(p,p_0) - s d(p_0,p) \\
        &\ge (1-s)r>0;
    \end{align*}
    hence $A\ge \phi(p_0) + (1-s)r>\phi(p_0)$. So, for $q \ne p_0$, 
    \[
    \phi(q)-d(q,p_0)\le \phi(p_0) < A \le \phi(q) + d(q,p_0); 
    \]
    so $|A-\phi(q)|\le d(q,p_0)$. 
    It follows that for $q \in M\setminus \{p_0\}$,
    \[
    \frac{|h(p_0)-h(q)|}{d(p_0,q)} = \frac{|A-\phi(q)|}{d(p_0,q)} \le 1.
    \]
    Thus, $h$ is $1$-Lipschitz. Since $\|\phi \restricted_{M\setminus \{p_0\}} \|_L=1$, we obtain that $\|h\|_L=1$ and $\nu_{L,h}(\phi)=1$. We have two subcases: 
    \begin{itemize}
        \item If the infimum in the definition of $A$ is attained at some $q$, then we can take $g:=h$ and $f:=\phi$. Indeed, $h$ is then strongly attaining its norm at $(p_0,q)$. Recall that $\phi$ does not strongly attain its norm. 
        \item If the infimum is not attained, then $h$ does not strongly attain its norm. Consider a sequence $(q_n)$ in $M\setminus \{p_0\}$ such that $\phi(q_n)+d(q_n,p_0)\to A$. Then 
            \[
            1\ge \frac{|h(p_0)-h(q_n)|}{d(p_0,q_n)} = \frac{|A-\phi(q_n)-d(q_n,p_0)+d(q_n,p_0)|}{d(p_0,q_n)} \to 1.
\]
Therefore, we can take $f=h$ and $g =d(\cdot ,p_0)-d(0,p_0)$.
    \end{itemize}
    In any of these two subcases, we find $f$ and $g$ that satisfy the three conditions in the statement.
\end{proof}

\begin{theorem}\label{thm:finite-metric-space-characterization}
Let $M$ be a pointed metric space. Then, the following assertions are equivalent:
\begin{enumerate}
    \item[\textup{(a)}] $M$ is finite.
    \item[\textup{(b)}] For every strongly norm-attaining $g\in S_{\Lip_0(M,\mathbb{K})}$,
    \[
    \LipNRA_g(M,\mathbb{K})=\Lip_0(M,\mathbb{K}).
    \]
    \item[\textup{(c)}] For every Banach space $Y$ and for every strongly norm-attaining $g\in S_{\Lip_0(M,Y)}$,
    \[
    \LipNRA_g(M,Y)=\Lip_0(M,Y).
    \]
    \item[\textup{(d)}] For every Banach space $Y$ and for every strongly norm-attaining $g\in S_{\Lip_0(M,Y)}$,
    \[\overline{\LipNRA_g(M,Y)}=\Lip_0(M,Y).
    \]
\end{enumerate}
\end{theorem}

\begin{proof}
(a) $\implies$ (c) follows from the finiteness of $M$ and the weak-star compactness of $B_{Y^*}$. (c) $\implies$ (b) is clear, and (b) $\implies$ (a) is a consequence of Proposition \ref{prop:non-isolated-point}. Since (c) $\implies$ (d) is also clear, it remains to prove (d) $\implies$ (a). Suppose that $M$ is infinite. As above, choose $\varphi\in S_{\Lip_0(M,\mathbb R)}$ that does not strongly attain its norm. Define $f,g:M\to \ell_\infty^2(\mathbb{K})$ by 
\[
g(p)=(\varphi(p), d(0,p))\quad\text{and}\quad f(p)=(\varphi(p),0)
\]
for all $p\in M$. Then $\|g\|_L=\|f\|_L=1$, and $\|g(p)-g(0)\|_\infty=d(0,p)$ for every $p\ne 0$. Choose distinct $p_n,q_n\in M$ such that 
\[
\frac{\varphi(p_n)-\varphi(q_n)}{d(p_n,q_n)} \to 1. 
\]
Then 
\[
e_1^* \left( \frac{g(p_n)-g(q_n)}{d(p_n,q_n)} \right) =e_1^* \left( \frac{f(p_n)-f(q_n)}{d(p_n,q_n)} \right)  =\frac{\varphi(p_n)-\varphi(q_n)}{d(p_n,q_n)} \to 1. 
\]
This implies that $\nu_{L,g}(f)=1$.

Let $(p,q,(a,b))\in\Pi_g$ with $(a,b)\in S_{\ell_1^2(\mathbb{K})}$ be given. Then 
\begin{align*}
1=\left|a \left(\frac{\varphi(p)-\varphi(q)}{d(p,q)}\right) +b \left(\frac{d(0,p)-d(0,q)}{d(p,q)}\right) \right|. 
\end{align*}
Since 
\[
\left|\frac{\varphi(p)-\varphi(q)}{d(p,q)}\right| <1, 
\]
we obtain $a=0$ and hence 
\[
\left\langle (a,b), \frac{f(p)-f(q)}{d(p,q)}\right \rangle =0.
\]Consequently, if $h\in \LipNRA_g(M,\ell_\infty^2(\mathbb{K}))$ attains its $g$-Lipschitz numerical radius at $(p,q,(a,b)) \in \Pi_g$, then
\[
\nu_{L,g}(h) \le \left|\frac{(h-f)(p)-(h-f)(q)}{d(p,q)}\right| 
\le\|h-f\|_L.
\]
It follows that
\[
1=\nu_{L,g}(f)\le\nu_{L,g}(h)+\|f-h\|_L\le2\|f-h\|_L.
\]
Therefore,
\[
\operatorname{dist}\left(f, \LipNRA_g(M,\ell_\infty^2(\mathbb{K}))\right)\ge\frac12>0;
\]
thus (d) fails. 
\end{proof}

\subsection{Length spaces and numerical-radius attainment}

Recall that a metric space $M$ is a \emph{length space} if, for every $p,q\in M$ and every $\eta>0$, there is a rectifiable curve joining $p$ and $q$ whose length is less than $d(p,q)+\eta$. The following theorem extends \cite[Theorem~2.4]{CJT}, which corresponds to $M=Y=X$ and $g=\id_X$ for a Banach space $X$.

\begin{theorem}\label{thm:non-density-length}
Let $M$ be a pointed length space, let $Y$ be a Banach space, and let $g\in S_{\Lipz(M,Y)}$. Then $\LipNRA_g (M,Y)$ is not dense in $\Lipz(M,Y)$.
\end{theorem}

\begin{proof}
We adapt the argument of Kadets, Mart\'in, and Soloviova \cite[Theorem 2.3]{KMS} to the present numerical-radius setting, using almost shortest curves as in \cite[Theorem 2.2]{CCGMR}. The key additional estimate controls the almost norming triples associated with the fixed reference map $g$. If $\Pi_g=\emptyset$, then $\LipNRA_g (M,Y)=\emptyset$. Suppose that $\Pi_g\ne\emptyset$ and fix $(p_0,q_0,y_0^*)\in\Pi_g$. Put
\[
u_0= \frac{g(p_0)-g(q_0)}{d(p_0,q_0)} \quad\text{and}\quad r=\Rea(y_0^*\circ g) \in \Lip_0 (M,\mathbb{R}).
\]
Then $\|u_0\|=\|r\|_L=1$, $y_0^*(u_0)=1$, and $r(p_0)-r(q_0)=d(p_0,q_0)$.
Choose a compact nowhere dense set $C\subset(r(q_0),r(p_0))$ with positive Lebesgue measure and define
\[
\Phi(t)=\int_0^t\chi_C(s)\,ds, 
\]
and also define $f\in\Lipz(M,Y)$ by
\[
f(p)=\Phi(r(p)) \, u_0
\]
for every $p\in M$. Then $\|f\|_L\le1$, and
\[
\vL{g}(f)\ge \left|y_0^* \left( \frac{f(p_0)-f(q_0)}{d(p_0,q_0)} \right) \right|=\frac{|C|}{d(p_0,q_0)}=:c>0.
\]

Suppose that $h\in\LipNRA_g (M,Y)$ satisfies $\rho:=\|h-f\|_L<\tfrac{c}{2}$. Then $\vL{g}(h)\ge c-\rho>\rho$. Choose $(p,q,y^*)\in\Pi_g$ at which $h$ attains its radius. Since
\[
|y^*(f(p)-f(q))|\ge(\vL{g}(h)-\rho)d(p,q)>0,
\]
we have $r(p)\ne r(q)$. Replacing $(p,q,y^*)$ by $(q,p,-y^*)$ if necessary, assume $r(q)<r(p)$. Choose $r(q)<a<b<r(p)$ with $[a,b]\cap C=\emptyset$.

For each $n$, choose a curve 
$\gamma_n:[0,1]\to M$ with $\gamma_n(0)=q$, $\gamma_n(1)=p$ and length less than $d(p,q)+\tfrac1n$.
By continuity of $r\circ\gamma_n$, choose $s_n < t_n\in(0,1)$ with \begin{equation}\label{eq:r(gamma_n(s_n)=a-and-r(gamma_n(t_n)=b}
    r(\gamma_n(s_n))=a  \quad \text{and} \quad  r(\gamma_n(t_n))=b.
\end{equation}
Also, we put 
\begin{equation}\label{eq:gamma_n(s_n)=u_n-gamma_n(t_n)=v_n}
    \gamma_n(s_n)=u_n \quad   \text{and}  \quad  \gamma_n(t_n)=v_n.
\end{equation}
By the choice of $\gamma_n$, we have
\[
d(q,u_n)+d(u_n,v_n)+d(v_n,p)<d(p,q)+\frac1n.
\]
Since $\|r\|_L=1$, $d(q,u_n)\ge a-r(q)>0$, $d(u_n,v_n)\ge b-a>0$ and $d(v_n,p)\ge r(p)-b>0$.

Using $y^*(g(p)-g(q))=d(p,q)$, we obtain
\begin{align*}
    &\left[d(q,u_n)-\Rea y^*(g(u_n)-g(q))\right]\\
&\quad+\left[d(u_n,v_n)-\Rea y^*(g(v_n)-g(u_n))\right]\\
&\quad+\left[d(v_n,p)-\Rea y^*(g(p)-g(v_n))\right]\\
&=d(q,u_n)+d(u_n,v_n)+d(v_n,p)-d(p,q)<\frac1n.
\end{align*}
In particular, we have

\begin{equation}\label{eq:est1}
    0\le1-\Rea y^*\left(\frac{g(u_n)-g(q)}{d(u_n,q)}\right)<\frac{1}{n\,d(u_n,q)}\le\frac{1}{n(a-r(q))},
\end{equation}
\begin{equation}\label{eq:est2}
    0\le1-\Rea y^*\left(\frac{g(p)-g(v_n)}{d(p,v_n)}\right)<\frac{1}{n\,d(p,v_n)}\le\frac{1}{n(r(p)-b)}.
\end{equation}
Fix $0<\delta<1$. Then $(u_n,q,y^*),(p,v_n,y^*)\in\Pi_g(\delta)$ for all sufficiently large $n$. Hence,
\[
\begin{aligned}
    |y^*(h(u_n)-h(q))|&\le\vL{g}^\delta(h) d(q,u_n),\\
    |y^*(h(p)-h(v_n))|&\le\vL{g}^\delta(h) d(v_n,p).
\end{aligned}
\]
Furthermore, from $\Phi(b)-\Phi(a)=\int_a^b \chi_C(s)\,ds=0$ and \eqref{eq:r(gamma_n(s_n)=a-and-r(gamma_n(t_n)=b} with \eqref{eq:gamma_n(s_n)=u_n-gamma_n(t_n)=v_n}, we see that 
\[
f(u_n)=f(\gamma_n(s_n))=\Phi(a) u_0 = \Phi(b) u_0 = f(\gamma_n(t_n))=f(v_n).
\]
So,
\begin{equation}\label{eq:est3}
    |y^*(h(v_n)-h(u_n))|=|y^*((h-f)(v_n)-(h-f)(u_n))|\le\rho d(u_n,v_n).
\end{equation}
By \eqref{eq:est1}, \eqref{eq:est2} and \eqref{eq:est3}, we have
\begin{align*}
    \vL{g}(h) &d(p,q) \\
    &=|y^*(h(p)-h(q))|\\
&\le|y^*(h(p)-h(v_n))|+|y^*(h(v_n)-h(u_n))|+|y^*(h(u_n)-h(q))|\\
&\le\vL{g}^\delta(h)\bigl(d(q,u_n)+d(v_n,p)\bigr)+\rho d(u_n,v_n)\\
&\le\vL{g}^\delta(h)\left(d(p,q)+\frac1n\right)-(\vL{g}^\delta(h)-\rho)(b-a),
\end{align*}
since $\vL{g}^\delta(h)\ge\vL{g}(h)>\rho$ and $d(u_n,v_n)\ge b-a$.
Consequently, letting $n\to\infty$ first and then $\delta \downarrow 0$, we have
\[
\vL{g}(h)d(p,q)\le\vL{g}(h)d(p,q)-(\vL{g}(h)-\rho)(b-a)<\vL{g}(h)d(p,q).
\]
This yields a contradiction.
\end{proof}

Recall that a pair of distinct points $(u,v)$ has \emph{property $(Z)$} if, for every $\varepsilon>0$, there exists $z\in M\setminus\{u,v\}$ such that
\[
d(u,z)+d(z,v)\le d(u,v)+\varepsilon\min\{d(u,z),d(v,z)\}.
\]
The space $M$ has property $(Z)$ if every pair of distinct points has this property; see \cite[Definition~1.1]{AMC}. It is also known in \cite[Theorem~5.4]{GLPRZ} that $(u,v)$ does not have the property $(Z)$ if and only if the molecule $m_{u,v}$ is a strongly exposed point of $B_{\mathcal{F}(M)}$. Finally, a complete metric space $M$ is length if and only if it has property $(Z)$ \cite[Main Theorem]{AMC}.  

For complete metric spaces, the preceding obstruction is exact. If $M$ is not a length space, a single reference map can be chosen so that every Lipschitz map attains its numerical radius. 

\begin{theorem}\label{thm:nonlength-lnra}
Let $M$ be a complete pointed metric space and $Y$ be a Banach space. The following assertions are equivalent:
\begin{enumerate}
 \item[\textup{(a)}] $M$ is a length space.
 \item[\textup{(b)}] $\overline{\LipNRA_g (M,Y)} \neq \Lip_0(M,Y)$ for every $g\in S_{\Lipz(M,Y)}$. 
\end{enumerate}
Moreover, if $M$ is not length, then there exists $g\in S_{\Lipz(M,Y)}$ such that $\LipNRA_g (M,Y)=\Lipz(M,Y)$.
\end{theorem}

For complete $M$, the length-space property is also equivalent to the Daugavet property of the real Lipschitz-free space $\mathcal{F}(M)$ and of $\Lip_0(M,\mathbb{R})$ \cite[Theorem~3.5]{GLPRZ}. Theorem \ref{thm:nonlength-lnra} therefore identifies the same metric geometry through numerical-radius attainment. 

\begin{proof}[Proof of Theorem \ref{thm:nonlength-lnra}]
The implication (a) $\implies$ (b) follows from Theorem~\ref{thm:non-density-length}.

It suffices to prove the ``moreover'' part. Suppose that $M$ is not a length space. Then  some pair $(u,v)$ fails property~$(Z)$ by \cite[Main Theorem]{AMC}.
By \cite{GLPRZ}, there is a function $\varphi\in S_{\Lip_0(M,\mathbb R)}$ \emph{peaking} at $(u,v)$. Thus, $\varphi(u)-\varphi(v)=d(u,v)$, and for every sequence of pairs $(p_n,q_n)$ with $p_n\ne q_n$, we have
\[
\frac{\varphi(p_n)-\varphi(q_n)}{d(p_n,q_n)}\to 1
\quad\implies\quad p_n\to u,\quad\text{and}\quad q_n\to v.
\]
Fix $y_0\in S_Y$ and put $g(x)=\varphi(x)y_0$ for $x\in M$.

Since $\|\varphi\|_L=1$ and $\|y_0\|=1$, we have $\|g\|_L=1$.
Fix $f\in\Lip_0(M,Y)$. By the definition of $\nu_{L,g}(f)$, we may
choose $(p_n,q_n,y_n^*)\in\Pi_g\left(\tfrac1n\right)$ for $n\ge 2$ such that
\[
\left|y_n^*\left(\frac{f(p_n)-f(q_n)}{d(p_n,q_n)}\right)\right|
\to \nu_{L,g}(f).
\]
Replacing $(p_n,q_n,y_n^*)$ by $(q_n,p_n,-y_n^*)$ whenever
necessary, we may assume that $\varphi(p_n)\ge\varphi(q_n)$
for every $n$. Put
\[
a_n=\frac{\varphi(p_n)-\varphi(q_n)}{d(p_n,q_n)}\in[0,1]
\]
and observe that \[
1-\frac1n < a_n\Rea y_n^*(y_0)\le a_n\le 1.
\]
It follows that $a_n\to1$ and $\Rea y_n^*(y_0)\to 1$. As $|y_n^*(y_0)|\le 1$, we also have $y_n^*(y_0)\to 1$.

Since $\varphi$ peaks at $(u,v)$, the convergence $a_n\to1$ gives $p_n\to u$ and $q_n\to v$. Thus, by continuity of $f$ and the fact
that $u\ne v$,
\[
\frac{f(p_n)-f(q_n)}{d(p_n,q_n)} \to z:=\frac{f(u)-f(v)}{d(u,v)} \quad \text{ in }Y.
\]
Since $\|y_n^*\|=1$ for every $n$, it follows that $|y_n^*(z)|\to\nu_{L,g}(f)$.

By weak-star compactness of $B_{Y^*}$, the sequence $(y_n^*)$ has
a subnet converging in the weak-star topology to some $y^*\in B_{Y^*}$.
Passing to the limit along this subnet gives
\[
y^*(y_0)=1, \quad\text{and}\quad |y^*(z)|=\nu_{L,g}(f).
\]
In particular, $\|y^*\|=1$. Since $\varphi(u)-\varphi(v)=d(u,v)$,
we have
\[
y^*\left(\frac{g(u)-g(v)}{d(u,v)}\right)=y^*(y_0)=1.
\]
Thus $(u,v,y^*)\in\Pi_g$, and $f$ attains its $g$-Lipschitz numerical radius at this triple. As $f$ was arbitrary, $\LipNRA_g(M,Y)=\Lip_0(M,Y)$.
\end{proof}

\section{Geometric estimates for the Lipschitz numerical index}

\subsection{Lipschitz maps with uniformly smooth range}\label{subsec:unif-smooth-range}

In this subsection, we obtain upper bounds for the $g$-Lipschitz numerical index by composing $g$ with suitable real-linear operators on a uniformly smooth target space. This yields an estimate for $\ell_p(\Lambda)$ and, in particular, shows that the index vanishes for Hilbert-valued maps whenever the target space has dimension at least two.

Throughout this subsection, if $Y$ is a \emph{smooth} Banach space, then for $u\in S_Y$, we denote by $\pi_Y(u)$ the unique element in $S_{Y^*}$ that satisfies $\pi_Y(u)(u)=1$.

Recall that the \emph{modulus of smoothness} of a Banach space $Y$ is defined by
\[
\rho_Y(t)=\sup_{\substack{x\in S_Y\\y\in B_Y}}\left(\frac{\|x+ty\|+\|x-ty\|}{2}-1\right)
\]
for $t\ge0$. The space $Y$ is \emph{uniformly smooth} if $\tfrac{\rho_Y(t)}{t}\to0$ as $t\downarrow0$.

We include a short proof of the following standard fact for completeness.
\begin{lemma}\label{lem:uniform-norm-expansion}
Let $Y$ be a real or complex uniformly smooth Banach space. Then
\[
\|u+\alpha v\|=1+\alpha\Rea\pi_Y(u)(v)+o(|\alpha|)
\]
as the real parameter $\alpha$ tends to zero, uniformly for $u\in S_Y$ and $v\in B_Y$.
\end{lemma}
\begin{proof}
Fix $u\in S_Y$, $v\in B_Y$, and $\alpha\in\mathbb R$. We have
\[
\|u\pm\alpha v\|\ge\Rea\pi_Y(u)(u\pm\alpha v)=1\pm\alpha\Rea\pi_Y(u)(v).
\]
On the other hand, the definition of $\rho_Y$ gives
\[
\|u+\alpha v\|+\|u-\alpha v\|\le2+2\rho_Y(|\alpha|).
\]
Combining these inequalities, we obtain
\[
\|u+\alpha v\|\le 2+2\rho_Y(|\alpha|)-\|u-\alpha v\|\le 1+\alpha\Rea\pi_Y(u)(v)+2\rho_Y(|\alpha|).
\]
Together with the lower bound for $\|u+\alpha v\|$, this gives
\[
0\le\|u+\alpha v\|-1-\alpha\Rea\pi_Y(u)(v)\le2\rho_Y(|\alpha|).
\]
Consequently,
\[
\sup_{\substack{u\in S_Y\\v\in B_Y}}\frac{\left|\|u+\alpha v\|-1-\alpha\Rea\pi_Y(u)(v)\right|}{|\alpha|}\le\frac{2\rho_Y(|\alpha|)}{|\alpha|}\to 0,
\]
which proves the asserted uniform expansion.
\end{proof}

\begin{proposition}\label{prop:lp-isometric-composition}
Let $Y$ be a real or complex uniformly smooth Banach space,
and let $F:Y\to Y$ be a bounded $\mathbb R$-linear operator.
Suppose that there exists $C\ge0$ such that
\[
|\pi_Y(u)(Fu)|\le C\|Fu\|
\]
for every $u\in S_Y$. Then, for every pointed metric space $M$ and every
$g\in S_{\Lipz(M,Y)}$, we have
\[
\nu_{L,g}(F\circ g)\le C\|F\circ g\|_L.
\]
In particular, if $F\circ g\ne0$, then $n_{L,g}(M,Y)\le C$.
\end{proposition}

\begin{proof}
Put $f=F\circ g$.
For $z\in B_Y\setminus\{0\}$, write $u=\tfrac{z}{\|z\|}$.
The estimate in the proof of
Lemma~\ref{lem:uniform-norm-expansion} gives, for every
$\alpha>0$ and $\theta\in\T$,
\[
\begin{aligned}
\|z+\alpha\theta Fz\|
&=\|z\|\|u+\alpha\theta Fu\|\\
&\le \|z\|
 +\alpha\|z\|\Rea\bigl(\theta\pi_Y(u)(Fu)\bigr)
 +2\|z\|\rho_Y(\alpha\|F\|)\\
&\le \|z\|+\alpha C\|Fz\|
 +2\rho_Y(\alpha\|F\|).
\end{aligned}
\]
The same inequality holds trivially when $z=0$.
Applying this estimate to
\[
z=\frac{g(p)-g(q)}{d(p,q)}
\]
for $p,q\in M$ with $p\neq q$
and using the real-linearity of $F$, we obtain
\[
\|g+\alpha\theta f\|_L
\le 1+\alpha C\|f\|_L+2\rho_Y(\alpha\|F\|).
\]
Since $Y$ is uniformly smooth,
$2\rho_Y(\alpha\|F\|)=o(\alpha)$ as $\alpha\downarrow0$.
Thus Corollary~\ref{cor:norm-derivative} yields
\[
\nu_{L,g}(f)
=\lim_{\alpha \to 0^+}
 \max_{\theta\in\T}
 \frac{\|g+\alpha\theta f\|_L-1}{\alpha}
\le C\|f\|_L. \qedhere
\]
\end{proof}

Real-linearity is sufficient because $F \circ g$ is used as a Lipschitz perturbation. In the complex case, this allows the conjugate-linear two-coordinate map used below and explains why the Hilbert-valued conclusion differs from the classical complex numerical index. 

We apply Proposition~\ref{prop:lp-isometric-composition} to $\ell_p(\Lambda)$ to prove part~(a) of Theorem~\ref{thm:intro-geometry}.

\begin{corollary}\label{cor:lp-relative-index-bound}
Let $1<p<\infty$, let $\Lambda$ be a set with $|\Lambda|\ge2$,
and let $\ell_p(\Lambda)$ be real or complex.
For every pointed metric space $M$ and every
$g\in S_{\Lipz(M,\ell_p(\Lambda))}$, we have
\[
n_{L,g}(M,\ell_p(\Lambda))\le M_p,\quad\text{and}\quad
M_p:=\max_{0\le t\le1}
\frac{|t^{p-1}-t|}{1+t^p}.
\]
In particular, if $H$ is a real or complex Hilbert space with $\dim H \geq 2$, then for every $g\in S_{\Lipz(M,H)}$, we have
\[
n_{T_g}(\F_{\mathbb{K}} (M),H)=n_{L,g}(M,H)=0.
\]
\end{corollary}

\begin{proof}
Write
$g(x)=(g_\lambda(x))_{\lambda\in\Lambda}$ for $x\in M$. Since $\|g\|_L=1$, there exist $p_0,q_0\in M$ and $i\in\Lambda$ such that $g_i(p_0)\ne g_i(q_0)$. Fix $j\in\Lambda\setminus\{i\}$ and define $F:\ell_p(\Lambda)\to \ell_p(\Lambda)$ by
\[
(Fx)_k=
\begin{cases}
\overline{x_j}, & k=i,\\
-\overline{x_i}, & k=j,\\
0, & k\notin\{i,j\}.
\end{cases}
\]
Here complex conjugation is understood as the identity in the real case. Then $F$ is a real-linear contraction. Moreover, $f:=F\circ g$ satisfies
\[
f_j(p_0)-f_j(q_0)
=-\overline{g_i(p_0)-g_i(q_0)}\ne0,
\]
so $0<\|f\|_L\le1$.

For $x\in S_{\ell_p(\Lambda)}$, put $a=|x_i|$ and $b=|x_j|$. Using the norming functional of $\ell_p(\Lambda)$ and the definition of $M_p$, we obtain
\[
\begin{aligned}
|\pi_{\ell_p(\Lambda)}(x)(Fx)|
=|a^{p-1}b-ab^{p-1}|\le M_p(a^p+b^p)=M_p\|Fx\|^p\le M_p\|Fx\|.
\end{aligned}
\]
Since $\ell_p(\Lambda)$ is uniformly smooth, Proposition~\ref{prop:lp-isometric-composition} gives
\[
n_{L,g}(M,\ell_p(\Lambda))
\le \nu_{L,g}\left(\frac{f}{\|f\|_L}\right)
\le M_p.
\qedhere
\]
\end{proof}

\subsection{Banach spaces with property $\beta$}\label{subsec:beta}

In this subsection, we establish lower bounds for the Lipschitz numerical index of Banach spaces with property $\beta$. Our results extend the corresponding estimates for the classical numerical index in \cite[Propositions 5 and 6]{FMP2003} to the Lipschitz setting and answer the question posed in \cite[Question 4.11]{CJT}. 

Recall that a Banach space $X$ has the \emph{property $\beta$ with constant $\rho$} ($0\le \rho <1$) if there is a family 
\[
\{(x_i,x_i^*):i\in I\} \subset X \times X^*
\]
satisfying 
\begin{enumerate}[(i)]
    \item $\|x_i\|=\|x_i^*\|= x_i^*(x_i)=1$ \text{ for } $i \in I$,
    \item $|x_i^* (x_j)| \leq \rho$ for $i,j \in I$ with $i\neq j$, 
    \item $\|x\| = \sup\{|x_i^*(x)|:i \in I\}$ for $x \in X$.
\end{enumerate}

The following proposition proves part~(b) of Theorem~\ref{thm:intro-geometry}.

\begin{proposition}\label{prop:property-beta}
    Let $X$ be a Banach space that has property $\beta$ with constant $\rho$. 
    \begin{enumerate}
        \item[\textup{(a)}] If $X$ is a real Banach space, then     \[
    n_L(X) \geq \frac{1-\rho}{1+\rho}.
    \]
    \item[\textup{(b)}] If $X$ is a complex Banach space, then 
    \[
    n_L(X) \geq 1-\rho.
    \]
    \end{enumerate}

\end{proposition}

\begin{proof}
We adapt the constructions in \cite[Propositions 5 and 6]{FMP2003} by realizing them as two consecutive increments of a Lipschitz map. 
    Let $\{(x_i,x_i^*) : i \in I\}$ be the family appearing in the definition of property $\beta$. 
    
    (a) Fix $\varepsilon>0$ and $f \in \Lip_0 (X,X)$, and take $p\neq q \in X$ and $j \in I$ such that 
    \begin{equation}\label{eq:choice-of-p-q-and-j}
    \left| x_j^* \left( \frac{f(p)-f(q)}{\|p-q\|} \right) \right| > \|f\|_L-\varepsilon.
    \end{equation}
    Let us write
    \begin{gather*}
    x=\frac{p-q}{\|p-q\|}\in S_X,\quad
    \lambda=\frac{1-\rho}{1+\rho}\in(0,1],\quad
    t=\frac{1+\lambda x_j^*(x)}{2}\in[0,1],\\
    z=\lambda x+(1-\lambda x_j^*(x))x_j,
    \quad\text{and}\quad
    w=\lambda x-(1+\lambda x_j^*(x))x_j.
    \end{gather*}
    
    We clearly have $x_j^*(z)=1$, $x_j^*(w)=-1$, and $|x_i^*(z)|\leq 1$, $|x_i^*(w)|\leq 1$ for $i\neq j$. This shows that $\|z\|=\|w\|=1$, with $x_j^*$ attaining its norm at both points. Moreover, note that 
    \begin{equation}\label{eq:convex-combination-of-z-w}
      \lambda x = tz+(1-t)w.
    \end{equation}
    Let us consider 
    \[
    r:= q+ \frac{\|p-q\|}{\lambda} tz.
    \]
    Then we have that 
    \begin{align*}
        p-r=(p-q)+(q-r)&=\|p-q\|x - \frac{\|p-q\|}{\lambda}tz  \\
        &=\|p-q\|x - \frac{\|p-q\|}{\lambda} (\lambda x -(1-t)w) \\
        &=\frac{\|p-q\|}{\lambda }(1-t)w.
    \end{align*}
    Moreover, 
    \[
    \|r-q\|=x_j^*(r-q), \quad\text{and} \quad \|p-r\|=(-x_j^*)(p-r).
    \]
    Observe that 
    \[
    |x_j^*(f(r)-f(q))| \le \nu_L(f) \|r-q\|, \quad \text{and} \quad |x_j^*(f(p)-f(r))| \le \nu_{L} (f) \|p-r\|. 
    \]
    Since 
    \[
    \|p-q\|=t\|p-q\|+(1-t)\|p-q\|=\lambda\|r-q\|+\lambda\|p-r\|,
    \]
    we have that 
    \begin{align*}
|x_j^*(f(p)-f(q))| &\le |x_j^*(f(p)-f(r)) | + |x_j^* (f(r)-f(q))| \\
&\le \nu_L(f) (\|p-r\|+\|r-q\|) = \frac{\nu_L (f)}{\lambda}\|p-q\|.
    \end{align*}
    By \eqref{eq:choice-of-p-q-and-j}, we then have 
    \[
    \lambda (\|f\|_L-\varepsilon) \le \nu_L (f). 
    \]
    Letting $\varepsilon\downarrow0$, we obtain that $\nu_L(f) \geq \lambda \|f\|_L$.

    (b) We fix $f \in \Lip_0(X,X)$, $p\neq q \in X$ and $j \in I$ as in \eqref{eq:choice-of-p-q-and-j}. Write 
    \[
    x=\frac{p-q}{\|p-q\|} \in S_X, \quad \text{and} \quad  \lambda =1-\rho\in(0,1]. 
    \]
    Considering, if necessary, $\theta \in \mathbb{T}$ such that 
    \[
    (\theta x_j, \overline{\theta} x_j^*) 
    \]
    satisfies that $\overline{\theta}x_j^*(x) \in [0,1]$, we may assume that 
    $$t:=x_j^*(x)\in[0,1].$$
    Let 
    \[
    a = \sqrt{1-\lambda^2t^2}, \quad z = \lambda x + i ax_j,\quad \text{and}\quad w= \lambda x - i a x_j.
    \]
    Then $x_j^*(z)=\lambda t+ia$ and $x_j^*(w)=\lambda t-ia$; so 
    \[
    |x_j^*(z)|=|x_j^*(w)|=1.
    \]
    For $i\neq j$, we also have 
    \[
    |x_i^*(z)|\leq \lambda |x_i^* (x)|+ a|x_i^*(x_j)|\leq \lambda + a\rho \leq 1.
    \]
    Similarly, $|x_i^*(w)|\leq 1$; thus, $\|z\|=\|w\|=1=|x_j^*(z)|=|x_j^*(w)|$, and 
    \begin{equation}
        z+w=2\lambda x.
    \end{equation}
    Consider 
    \[
    r:=q+\frac{\|p-q\|}{2\lambda } z.
    \]
    Arguing as in part (a), we observe that 
    \begin{itemize}
        \item $\left |x_j^* \left( f(r)-f(q) \right) \right| \leq \nu_L(f) \|r-q\|$; 
        \item $\left |x_j^* \left( {f(p)-f(r)} \right) \right| \leq \nu_L(f)\|p-r\|$;
        \item $\|r-q\|+\|p-r\| =\frac{\|p-q\|}{\lambda}.$
    \end{itemize}     

    It follows that 
    \begin{align*}
\lambda(\|f\|_L-\varepsilon)
    &< \frac{\lambda}{\|p-q\|}\bigl|x_j^*(f(p)-f(q))\bigr| \\
    &\le \frac{\lambda\nu_L(f)}{\|p-q\|} \bigl(\|p-r\|+\|r-q\|\bigr)=\nu_L(f).
    \end{align*}
    Letting $\varepsilon\downarrow0$, we obtain $\nu_L(f) \geq \lambda \|f\|_L$.
\end{proof}

It is natural to ask whether the lower bounds in Proposition~\ref{prop:property-beta} remain valid for $n_{L,G}(X,X)$ when the identity is replaced by an arbitrary norm-one operator $G$. The following result shows that the index may vanish even when $G$ is arbitrarily close to the identity.

\begin{proposition}\label{prop:property-beta-instability}
    Let $X$ be a Banach space with $\dim (X)\geq 2$ and property $\beta$. For every $\varepsilon>0$, there exists a norm-one operator $G\in\mathcal{L}(X,X)$ such that $\|G-\id_X \| < \varepsilon$ and $n_{L,G}(X,X)=n_{G}(X,X)=0$.
\end{proposition}

\begin{proof}
        Let $\{(x_i,x_i^*) : i \in I\}$ be the family appearing in the definition of property $\beta$ (with constant $\rho$). Fix $j \in I$ and consider the rank-one projection $P$ on $X$ defined as 
        \[
        Px=x_j^*(x) x_j.
        \]
        Since $\dim(X)\geq 2$, there exists $v \in S_X \cap \ker x_j^*$. For $0<t<1$, define $T, G_t \in \mathcal{L}(X,X)$ by 
        \[
        Tx = x_j^*(x)v\quad\text{for every}~x\in X, \quad\text{and}\quad G_t = (1-t)\id_X + tP.
        \]
        Note that $\|T\|=\|G_t\|=1$ and $G_t - \id_X = (-t)(\id_X -P).$ So, $\|G_t-\id_X\| =t\|\id_X-P\|\to0$ as $t \to 0$. 

        Let $0<s \le t(1-\rho)$ and $\theta \in \mathbb{T}$ be given. We claim that $\|G_t+\theta sT\| =1.$ Let $x \in B_X$. Since $x_j^*(G_t x)=x_j^*(x)$ and $x_j^*(Tx) =0$, we have 
        \[
        x_j^* ((G_t + \theta s T)x)= x_j^*(G_tx)=x_j^*(x).
        \]
        For $i\neq j$, we have 
        \begin{align*}
            |x_i^* ((G_t + \theta s T)x)| &= |(1-t)x_i^*(x)+tx_j^*(x)x_i^*(x_j)+\theta s x_j^*(x)x_i^*(v)| \\
            &\leq (1-t)+t\rho + s \leq 1.
        \end{align*}
        It follows that $\|G_t+\theta s T\|=1$. By a characterization of the numerical index with respect to an operator \cite[Proposition 3.3]{KMMPQ}, we have 
        \[
        1=\max_{\theta\in\mathbb{T}} \|G_t + \theta sT\| \geq 1+ n_{G_t}(X,X)s,
        \]
        which shows $n_{G_t}(X,X)=0.$ Since $0\leq n_{L,G_t}(X,X)\leq n_{G_t}(X,X)=0$, we complete the proof.
\end{proof}

\section{Linear reference operators: equality and duality}\label{sec:comparison}
Throughout this section, all Banach spaces are real. We first prove part~(b) of Theorem~\ref{thm:intro-comparison}, treating separable domains in Theorem~\ref{thm:separable-domain} and targets that are $1$-complemented in their biduals in Theorem~\ref{thm:dual-target}. These proofs adapt the arguments of P\'erez-Hern\'andez \cite[Theorems~2.1 and~3.1]{PH} for the identity operator on separable spaces and on dual spaces, respectively. We then apply the dual-target case to the composition of $G$ with the canonical embedding $J_Y$ to obtain the lower bound in part~(a); the upper bound follows from Remark~\ref{rem:linear}(b).

\subsection{Separable domains}

\begin{lemma}\label{lem:differentiable-estimate}
Let $X$ be a separable Banach space, let $Y$ be a Banach space, and let $\gamma$ be a nondegenerate Gaussian probability measure on $X$ with mean zero such that $\int_X\|u\|\,d\gamma(u)<\infty$; such a measure exists by \cite[Lemma~2.3]{PH}. For $f\in\Lipz(X,Y)$, the map $S_\gamma f\colon X\to Y$ defined by
$$S_\gamma f(x)=\int_X\bigl[f(x+u)-f(u)\bigr]\,d\gamma(u)\quad\text{for every}~x\in X$$
belongs to $\Lipz(X,Y)$, satisfies $\|S_\gamma f\|_L\le\|f\|_L$, and is weakly G\^ateaux differentiable at every point of $X$. Moreover, $\nu_{L,G}^{\delta}(S_\gamma f)\le\nu_{L,G}^{\delta}(f)$ for every $G\in S_{\Lin(X,Y)}$ and every $0<\delta<1$.
\end{lemma}
\begin{proof}
For each $x\in X$, the map $u\mapsto f(x+u)-f(u)$ is $2\|f\|_L$-Lipschitz and bounded by $\|f\|_L\|x\|$. Being continuous on the separable space $X$, it is separably valued and weakly measurable, hence strongly measurable by the Pettis measurability theorem \cite[Theorem~II.1.2]{DU}; consequently, it is Bochner integrable with respect to the probability measure $\gamma$ \cite[Theorem~II.2.2]{DU}. 

The remaining assertions except the last one follow as in \cite[Proposition~2.4]{PH}, since \cite[p.~129]{GK2003} applies to Lipschitz maps from a separable Banach space into any Banach space. 
To prove the numerical-radius estimate, fix $G\in S_{\Lin(X,Y)}$
and $0<\delta<1$. By linearity of $G$, for all $p,q,x_0\in X$
with $p\ne q$ and every $t>0$, we have
\begin{equation}\label{eq:invariance}
\langle G,m_{x_0+tp,x_0+tq}\rangle
=\frac{G(p-q)}{\|p-q\|}
=\langle G,m_{p,q}\rangle.
\end{equation}
Taking $t=1$, we see that if $(p,q,y^*)\in\Pi_G(\delta)$, then
$(p+u,q+u,y^*)\in\Pi_G(\delta)$ for every $u\in X$. Hence
\[
\bigl|y^*(\langle S_\gamma f,m_{p,q}\rangle)\bigr|
\leq \int_X \bigl|y^*(\langle f,m_{p+u,q+u}\rangle)\bigr|
\,d\gamma(u)
\leq \nu_{L,G}^{\delta}(f).
\]
Taking the supremum over $(p,q,y^*)\in\Pi_G(\delta)$, we obtain $\vL{G}^\delta(S_\gamma f)\le\vL{G}^\delta(f)$.
\end{proof}

\begin{theorem}\label{thm:separable-domain}
Let $X$ be a separable Banach space, $Y$ be a Banach space, and let $G \in \mathcal{L}(X,Y)$ be a norm-one operator. Then 
\[
n_{L,G} (X,Y) = n_G(X,Y).
\]
\end{theorem}
\begin{proof}
By Remark~\ref{rem:linear}(b), it suffices to prove
$n_G(X,Y)\le n_{L,G}(X,Y)$.
Fix $f\in\Lip_0(X,Y)$ with $\|f\|_L=1$ and
$0<\varepsilon<1$, and let $\gamma$ be as in
Lemma~\ref{lem:differentiable-estimate}.

By Lemma~\ref{lem:differentiable-estimate} and the argument in \cite[proof of Theorem~2.1]{PH}, there exist $N\in\mathbb N$ and $x_0\in X$ such that the map
\[
F_N(x):=\int_X\left[f\left(x+\frac{u}{N}\right)-f\left(\frac{u}{N}\right)\right]\,d\gamma(u)\quad\text{for every}~ x\in X,
\]
is weakly G\^ateaux differentiable everywhere and satisfies
\[
\|dF_N(x_0)\|>1-\varepsilon,\quad\text{and}\quad\nu_{L,G}(F_N)\le\nu_{L,G}(f).
\]
Since $dG(x_0)=G$, Proposition~\ref{prop:weak-derivative-formulas}(a) gives
\begin{align*}
n_G(X,Y)\le \frac{\nu_G(dF_N(x_0))}{\|dF_N(x_0)\|}\le \frac{\nu_{L,G}(F_N)}{\|dF_N(x_0)\|}\le \frac{\nu_{L,G}(f)}{1-\varepsilon}.
\end{align*}
Letting $\varepsilon\downarrow0$ and then taking the infimum over all such $f$, we obtain
\[
n_G(X,Y)\le n_{L,G}(X,Y).\qedhere
\]
\end{proof}

\begin{remark} 
    Theorem \ref{thm:separable-domain} fails in general for a Lipschitz function $g\in S_{\Lip_0(X,Y)}$. One may ask if the following natural extension/formulation of Theorem \ref{thm:separable-domain} holds: 
    \[
    n_{L,g} (X,Y) = n_{\frac{dg(u)}{\|dg(u)\|}} (X,Y)
    \]
    for $g\in S_{\Lip_0(X,Y)}$ with non-zero $dg(u) \in \mathcal{L}(X,Y)$ at $u \in X$. Even when $X=Y=\mathbb{R}$, we take $g(t) := \log (\frac{1+e^t}{2})$  and $f(t):= t-g(t)$ for $t \in \mathbb{R}$. Note that $g \in \Lip_0 (\mathbb{R},\mathbb{R})$ with 
    \[
    dg(u)=\frac{e^u}{1+e^u} \text{Id}_\R
    \]
    for every $u\in \R$.
    Noting that $f(t)+g(t)=t$, $g'(t)\to1$ as $t\to\infty$, and $f'(t)\to0$ as $t \to \infty$, we observe that $\nu_{L,g}(f)=0.$ Since $\|f\|_L=\|g\|_L=1$, we obtain $n_{L,g}(\mathbb{R},\mathbb{R})=0$ (or, simply use Corollary \ref{cor:interval} as $|g'(t)|<1$ for every $t \in \mathbb{R}$). On the other hand, 
    \[
    n_{\frac{dg(u)}{\|dg(u)\|}} (\mathbb{R},\mathbb{R}) = n_{\text{Id}_\mathbb{R}} (\mathbb{R},\mathbb{R})=1
    \]
    for every $u \in \mathbb{R}$.
\end{remark}

\subsection{Targets that are $1$-complemented in their biduals}

We now examine the role of the target space in the comparison between linear and Lipschitz numerical indices. Specifically, we consider Banach spaces $Y$ whose canonical images in $Y^{**}$ are $1$-complemented. This class includes dual Banach spaces, and the projection onto $Y$ allows constructions in the bidual to be brought back to the original target. We investigate the consequences of this structure for the associated numerical indices. In particular, the identity case extends the equality $n_L (Y) = n(Y)$ from real dual spaces to all real spaces that are $1$-complemented in their biduals.

The following two lemmas are straightforward adaptations of the predual reduction and invariant-mean arguments in \cite[Section~3]{PH}. Lemma~\ref{lem:predual} follows from Goldstine's theorem. For Lemma~\ref{cor:linearize-from-invariant-means}, one applies invariant means first over $W$ and then over $X$. Since $G$ is linear, the almost-norming condition defining $\Pi_G(\delta)$ is preserved under simultaneous translations of the two points. Consequently, the same averaging estimates hold for each $0<\delta<1$. We omit the proofs.

\begin{lemma}\label{lem:predual}
Let $X$ and $Y$ be Banach spaces, let $G\in S_{\Lin(X,Y^*)}$, and let $f\in\Lipz(X,Y^*)$. For $0<\delta<1$, 
\begin{equation}\label{eq:in-lem:predual}
\nu_{L,G}^{\delta}(f)=\sup\left\{\left|\left(\frac{f(p)-f(q)}{\|p-q\|}\right)(y)\right| : \begin{aligned} &p,q\in X,\ p\neq q,\ y\in S_Y, \\
     &\left(\frac{G(p)-G(q)}{\|p-q\|}\right)(y) >1-\delta
\end{aligned} \right\}.
\end{equation}
That is, to compute $\nu_{L,G}^\delta(f)$, it is enough to consider elements of the predual $Y$ instead of elements of $Y^{**}$.
\end{lemma}

\begin{lemma}\label{cor:linearize-from-invariant-means}
Let $X$ and $Y$ be Banach spaces, let $G\in S_{\Lin(X,Y^*)}$, let $W\subseteq X$ be a linear subspace, and let $f\in\Lipz(X,Y^*)$ be such that $f\restricted_W$ is linear. Then there exists $T\in\Lin(X,Y^*)$ such that, for every $0<\delta<1$,
\[
T\restricted_W=f\restricted_W,\quad\|T\|\le\|f\|_L,\quad\text{and}\quad\nu_G^\delta(T)\le\vL{G}^\delta(f).
\]
\end{lemma}

\begin{theorem}\label{thm:dual-target}
    Let $X$ and $Y$ be Banach spaces. Suppose that $Y$ is $1$-complemented in $Y^{**}$. Then for every $G \in S_{\mathcal{L}(X,Y)}$, we have 
    \[
    n_{L,G} (X,Y) = n_G (X, Y).
    \]
    More precisely, for every $f \in S_{\Lip_0(X, Y)}$ and $0<\varepsilon<1$, there exists $T \in \mathcal{L}(X,Y)$ such that 
    \[
    1-\varepsilon<\|T\|\leq 1, \quad \text{and} \quad \nu_G (T) \leq \nu_{L,G} (f).
    \]
\end{theorem}
\begin{proof}
By Remark~\ref{rem:linear}(b), it is enough to show $n_G(X,Y)\le n_{L,G}(X,Y)$.

Let $P\in\Lin(Y^{**},Y)$ satisfy $\|P\|=1$ and $P\circ J_Y=\id_Y$. Fix $f\in S_{\Lip_0(X,Y)}$ and
$0<\varepsilon<1$. Take a closed separable subspace $Z\subseteq X$ such that $\|f\restricted_Z\|_L=1$.
Following the construction in \cite[proof of Theorem~3.1]{PH}, with Lemma~\ref{lem:differentiable-estimate} in place of \cite[Proposition~2.4]{PH}, let $F_N\colon X\to Y$ be the smoothed map
\[
F_N(x)=\int_Z
\left[f\left(x+\frac{u}{N}\right)-f\left(\frac{u}{N}\right)\right]
\,d\gamma(u)
\]
for every $x\in X$, where $\gamma$ is the Gaussian measure on $Z$ and $N\in\mathbb N$ is chosen as in that construction. Let $z_0\in Z$ be the differentiation point, and put $D=d(F_N\restricted_Z)(z_0)\in\mathcal L(Z,Y)$.

Let $\mathcal U$ be an ultrafilter on $(0,\infty)$ containing $(0,a)$ for every $a>0$, and define $\psi\colon X\to Y^{**}$ by
\[
\langle\psi(x),y^*\rangle
=
\lim_{t,\mathcal U}
\frac{y^*\bigl(F_N(z_0+tx)-F_N(z_0)\bigr)}{t}
\]
for every $x\in X$ and $y^*\in Y^*$. Since $\|F_N\|_L\le1$, the above formula defines an element $\psi(x)\in Y^{**}$. As in the cited construction, $\psi\in\Lip_0(X,Y^{**})$ and
\[
\|D\|>1-\varepsilon,\quad\|\psi\|_L\le1,\quad\text{and}\quad\psi\restricted_Z=J_YD.
\]
Moreover, Lemma~\ref{lem:predual} and \eqref{eq:invariance} give $\nu_{L,J_Y\circ G}^{\delta}(\psi)\le\nu_{L,G}^{\delta}(f)$ for every $0<\delta<1$.

By Lemma~\ref{cor:linearize-from-invariant-means}, there exists $A\in\mathcal L(X,Y^{**})$ such that
\[
A\restricted_Z=J_YD,\quad\|A\|\le1,\quad\text{and}\quad\nu_{J_Y\circ G}^{\delta}(A)\le \nu_{L,J_Y\circ G}^{\delta}(\psi)
\]
for every $0<\delta<1$. Set $T=PA$. Since $T\restricted_Z=D$, we have $1-\varepsilon<\|D\|\le\|T\|\le1$.
Since $P^*$ is an isometry and $(P^*y^*)(J_Y\circ Gx)=y^*(Gx)$ for every $x\in X,\ y^*\in Y^*$, for every $0<\delta<1$, we have
\[
\nu_G^\delta(T)\le \nu_{J_Y\circ G}^\delta(A)\le \nu_{L,J_Y\circ G}^\delta(\psi)\le \nu_{L,G}^\delta(f).
\]
Letting $\delta\downarrow0$ gives $\nu_G(T)\le\nu_{L,G}(f)$, proving the stated operator estimate. In particular,
\[
n_G(X,Y)\le \frac{\nu_G(T)}{\|T\|}\le \frac{\nu_{L,G}(f)}{1-\varepsilon}.
\]
Letting $\varepsilon\downarrow0$ and taking the infimum over all such $f$ yields $n_G(X,Y)\le n_{L,G}(X,Y)$.
\end{proof}

\begin{corollary}\label{cor:unified}
    Let $X$ and $Y$ be Banach spaces, and let $G \in S_{\mathcal{L}(X,Y)}$. Suppose that either $X$ is separable, or $Y$ is $1$-complemented in its bidual. Then 
    \[
    n_{L,G}(X,Y)=n_G(X,Y).
    \]
    Consequently, the following are equivalent.
\begin{enumerate}[label=(\alph*)]
        \itemsep0.25em
\item[\textup{(a)}] $n_{L,G}(X,Y)=1$;
\item[\textup{(b)}] $n_G(X,Y)=1$; 
\item[\textup{(c)}] $G$ is a spear operator. 
    \end{enumerate}
\end{corollary}

\begin{proof}
    Apply Theorem \ref{thm:separable-domain} in the first case and Theorem \ref{thm:dual-target} in the second. The characterization of spear operators follows from \cite[Proposition 3.3]{KMMPQ}.
\end{proof}

\subsection{A duality bound for arbitrary real Banach spaces}

It remains to prove part~(a) of Theorem~\ref{thm:intro-comparison}. We first identify the numerical index with respect to $G^*$ with that of $J_Y\circ G$ and then prove the lower bound in Theorem~\ref{thm:duality-sandwich}.

The spaces $\Lin(X,Y^{**})$ and $\Lin(Y^*,X^*)$ are both canonically isometric to the space $\mathcal{B}(X\times Y^*)=(X\hat{\otimes}_\pi Y^*)^*$ of bounded bilinear forms on $X\times Y^*$ \cite[Section~2.2]{Ryan}. The resulting surjective linear isometry $\psi\colon\Lin(X,Y^{**})\to\Lin(Y^*,X^*)$ is given by $\psi(B)=B^*\circ J_{Y^*}$.

\begin{lemma}\label{lem:bidual-identification}
Let $X$ and $Y$ be Banach spaces, and let $G\in S_{\Lin(X,Y)}$. Then
\[
\nu_{G^*}(\psi(B))=\nu_{J_Y\circ G}(B)
\]
for every $B\in\Lin(X,Y^{**})$. Consequently, $n_{G^*}(Y^*,X^*)=n_{J_Y\circ G}(X,Y^{**})$.
\end{lemma}

\begin{proof}
Since $\psi$ is a surjective linear isometry and $J_Y\circ G\in S_{\Lin(X,Y^{**})}$, \cite[Lemma~2.4(b)]{KMMPQ} gives
\[
\nu\left(\Lin(Y^*,X^*),\psi(J_Y\circ G),\psi(B)\right)=\nu\left(\Lin(X,Y^{**}),J_Y\circ G,B\right) 
\]
for every $B\in\Lin(X,Y^{**})$. Since $\psi(J_Y\circ G)=G^*$, this is the first assertion. Since $\psi$ maps the unit sphere of $\Lin(X,Y^{**})$ onto that of $\Lin(Y^*,X^*)$, taking infima gives the second.
\end{proof}

\begin{theorem}\label{thm:duality-sandwich}
Let $X$ and $Y$ be Banach spaces, and let $G\in S_{\Lin(X,Y)}$. Then
\[
n_{G^*}(Y^*,X^*)\,\|f\|_L\leq\nu_{L,G}(f)
\]
for every $f\in\Lipz(X,Y)$. Consequently, $n_{G^*}(Y^*,X^*)\leq n_{L,G}(X,Y)\leq n_G(X,Y)$. In particular, $n(X^*)\leq n_L(X)\leq n(X)$ for every real Banach space $X$.
\end{theorem}

\begin{proof}
Let
\[
\mathcal J:\operatorname{Lip}_0(X,Y)
\longrightarrow \operatorname{Lip}_0(X,Y^{**}),
\quad\text{and}\quad
\mathcal Jf=J_Y\circ f.
\]
Since $\mathcal J$ is a linear isometry, the Hahn--Banach theorem and Proposition~\ref{prop:abstract-range} give
\[
\nu_{L,J_Y\circ G}(\mathcal Jf)=\nu_{L,G}(f)
\]
for every $f\in\Lipz(X,Y)$. Since $Y^{**}$ is a dual space, 
\[
\begin{aligned}
n_{G^*}(Y^*,X^*)\|f\|_L
&=n_{J_Y\circ G}(X,Y^{**})\|\mathcal Jf\|_L \quad (\text{by Lemma \ref{lem:bidual-identification}})\\
&=n_{L,J_Y\circ G}(X,Y^{**})\|\mathcal Jf\|_L \quad (\text{by Theorem \ref{thm:dual-target}})\\
&\leq \nu_{L,J_Y \circ G}(\mathcal Jf)
=\nu_{L,G}(f).
\end{aligned}
\]
Taking the infimum over all $f$ with $\|f\|_L=1$
and using Remark~\ref{rem:linear}(b) gives the comparison of the indices.
The final assertion follows by taking $Y=X$ and $G=\operatorname{Id}_X$.
\end{proof}

\begin{remark}
Recall that a Banach space $Y$ is \emph{$L$-embedded} if $$Y^{**}=J_Y(Y)\oplus_1 Y_s$$ 
for some closed subspace $Y_s\subseteq Y^{**}$. If $Y$ is real and $L$-embedded, then \cite[Proposition~6.17]{KMMPQ} yields
\[
n_{G^*}(Y^*,X^*) =n_G(X,Y)
\]
for every real Banach space $X$ and every $G\in S_{\mathcal L(X,Y)}$. Therefore, Theorem~\ref{thm:duality-sandwich} shows that 
\[
n_{G^*}(Y^*,X^*) =n_{L,G} (X,Y)=n_G(X,Y)
\]
for such $X, Y$ and $G$.

This applies to all real $L_1(\mu)$ spaces, which are
$L$-embedded \cite[Example~IV.1.1]{HWW}. In particular, this yields 
$$
n_{L,G}(X,L_1[0,1])=n_G(X,L_1[0,1])
$$ for every real Banach space $X$ and every norm-one operator
$G:X\to L_1[0,1]$. Since $L_1[0,1]$ is not isomorphic to a dual space \cite[Chapter~III]{DU}, choosing $X$ nonseparable provides examples beyond both the separable-domain and dual-target cases. Moreover, taking $Y=X$ and $G=\mathrm{Id}_X$, we obtain
\[
n_L(X)=n(X)=n(X^*)
\]
for every real $L$-embedded Banach space $X$. To the best of our knowledge, the equality $n_L(X)=n(X)$ is new in this generality.
\end{remark}

\section*{Acknowledgements} 

G. Choi was supported by the National Research Foundation of Korea (NRF) grant funded by the Korea government (MSIT) (RS-2026-25475373).
M. Han was supported by Global - Learning \& Academic research institution for Master’s·PhD students, and Postdocs (LAMP) Program of the National Research Foundation of Korea (NRF) grant funded by the Ministry of Education (RS-2024-00445180).
M. Jung was supported by the research fund of Hanyang University (HY-202500000003346) and by the National Research Foundation of Korea (NRF) grant funded by the Korean government (MSIT) (RS-2026-25599117).

\section*{Declaration of AI Use} 

During the development of this work, ChatGPT and Claude were used as supplementary research tools in Sections \ref{sec:metric-interpretation} and \ref{sec:graph} to explore graph-theoretic questions, examine examples, and discuss possible approaches. In addition, ChatGPT suggested an initial idea for the proof of the result concerning property $\beta$ in Section \ref{subsec:beta}. An earlier ChatGPT-assisted proof of Theorem~\ref{thm:duality-sandwich} was subsequently replaced by a different proof developed by the authors. Claude also provided a proof of the Hilbert space version of Proposition~\ref{prop:lp-isometric-composition}; based on its underlying idea, the authors generalized the result to uniformly smooth target spaces and derived Corollary~\ref{cor:lp-relative-index-bound}. All mathematical arguments presented in this paper were independently verified by the authors, who take full responsibility for their correctness. Apart from the uses described above, no AI tools were used to develop mathematical ideas or construct proofs in the remaining sections.

\end{document}